\documentclass[10pt]{article}
\usepackage{amssymb,amsmath,latexsym}
\begin{document}
\allowdisplaybreaks
\begin{doublespace}

\newtheorem{thm}{Theorem}[section]
\newtheorem{lemma}[thm]{Lemma}
\newtheorem{defn}[thm]{Definition}
\newtheorem{prop}[thm]{Proposition}
\newtheorem{corollary}[thm]{Corollary}
\newtheorem{remark}[thm]{Remark}
\newtheorem{example}[thm]{Example}
\numberwithin{equation}{section}
\def\ee{\varepsilon}
\def\qed{{\hfill $\Box$ \bigskip}}
\def\NN{{\cal N}}
\def\AA{{\cal A}}
\def\MM{{\cal M}}
\def\BB{{\cal B}}
\def\CC{{\cal C}}
\def\LL{{\cal L}}
\def\DD{{\cal D}}
\def\FF{{\cal F}}
\def\EE{{\cal E}}
\def\QQ{{\cal Q}}
\def\RR{{\mathbb R}}
\def\bR{{\mathbb R}}
\def\R{{\mathbb R}}
\def\L{{\bf L}}
\def\K{{\bf K}}
\def\S{{\bf S}}
\def\A{{\bf A}}
\def\E{{\mathbb E}}
\def\F{{\bf F}}
\def\P{{\mathbb P}}
\def\N{{\mathbb N}}
\def\eps{\varepsilon}
\def\wh{\widehat}
\def\wt{\widetilde}
\def\pf{\noindent{\bf Proof.} }
\def\beq{\begin{equation}}
\def\eeq{\end{equation}}
\def\bee{\begin{equation}}
\def\eee{\end{equation}}
\def\sA {{\cal A}}
\def\1{{\bf 1}}
\def\nn{\nonumber}

\title{\Large \bf
Two-sided Green function estimates for killed subordinate Brownian motions}

\author{{\bf Panki Kim}\thanks{Supported by Basic Science Research Program through
the National Research Foundation of Korea(NRF) grant funded by the
Korea government(MEST)(2010-0001984).} \quad {\bf Renming Song}
\quad and  \quad {\bf Zoran Vondra{\v{c}}ek}\thanks{Supported in
part by the MZOS grant 037-0372790-2801.} }

\date{February 8, 2011}

\maketitle

\begin{abstract}
A subordinate Brownian motion is a L\'evy process which can be
obtained by replacing the time of the Brownian motion by an
 independent subordinator. The infinitesimal generator
of a subordinate Brownian motion is $-\phi(-\Delta)$, where $\phi$
is the Laplace exponent of the subordinator. In this paper, we
consider a large class of subordinate Brownian motions without diffusion component
and with $\phi$ comparable to a regularly varying function at infinity.
This class of processes includes  symmetric stable
processes, relativistic stable processes, sums of independent
symmetric stable processes, sums of independent relativistic stable
processes, and much more. We give sharp two-sided estimates on the
Green functions of these subordinate Brownian motions in any bounded
$\kappa$-fat open set $D$. When $D$ is a bounded $C^{1,1}$ open set,
we establish an explicit form of the estimates in terms of the
distance to the boundary. As a consequence of such sharp Green
function estimates, we obtain a boundary Harnack principle
in $C^{1,1}$ open sets with explicit rate of decay.
\end{abstract}



\section{Introduction}
The investigation of fine potential-theoretic properties of
discontinuous Markov processes in the Euclidean space began in the
late 1990's with the study of symmetric stable processes. One of the
first results obtained in this area was sharp Green function
estimates of symmetric $\alpha$-stable processes in bounded
$C^{1,1}$ domains in $\R^d$, $0<\alpha<2$, $d\ge 2$. Recall that if
$X$ is a symmetric Markov process in $\R^d$ and $D$ is an open
subset of $\R^d$, then the Green function $G_D(x, y)$ of $X$ in $D$
(if it exists) is the density of the mean occupation time for $X$
before exiting $D$, that is, the density of the measure
$$
U\mapsto \E_x\int^{\tau_D}_0{\bf 1}_U(X_t)dt, \quad U\subset D,
$$
where $\tau_D$ is the first time the process $X$ exits $D$.
Analytically speaking, if $\mathcal L$ is the infinitesimal
generator of $X$ and $\mathcal L|_D$ is the restriction of $\mathcal
L$ to $D$ with zero exterior condition, then $G_D(\cdot, y)$ is the
solution of $(\mathcal L|_D) u=-\delta_y$.

A process $X=(X_t:\, t\ge 0)$ is called a (rotationally) symmetric
$\alpha$-stable (L\'evy) process, $0<\alpha<2$, if it is a L\'evy
process whose characteristic exponent $\Phi$, defined by $\E[\exp\{i
\theta \cdot X_t \}] =\exp\{-t\Phi(\theta)\}$, is given by
$\Phi(\theta)=|\theta|^{\alpha}$. The infinitesimal generator of a
symmetric $\alpha$-stable process is $-(-\Delta)^{\alpha/2}$. The
paths of the symmetric $\alpha$-stable process $X$ are purely
discontinuous, as opposed to the case $\alpha=2$ corresponding to
Brownian motion which has continuous paths. It was independently
shown in \cite{CS1} and \cite{Ku1}
that if $D$ is a bounded $C^{1,1}$ domain in $\R^d$, $G_D(x,y)$ the
Green function of the symmetric $\alpha$-stable process in  $D$, and
$\delta_D(x)$ the distance between the point $x$ and the complement
$D^c$ of $D$, then there exists a constant $c>1$ (depending on $D$
and $\alpha$) such that
\begin{equation}\label{e:green-for-stable}
c^{-1} \left(1\wedge
\frac{(\delta_D(x)\delta_D(y))^{\alpha/2}}{|x-y| ^\alpha}\right)
\frac{1}{|x-y|^{d-\alpha}}\le G_D(x,y) \le c \left(1\wedge
\frac{(\delta_D(x)\delta_D(y))^{\alpha/2}}{|x-y| ^\alpha}\right)
\frac{1}{|x-y|^{d-\alpha}}\, ,
\end{equation}
for all $x,y\in D$. Here and in the sequel, for $a, b\in \bR$,
$a\wedge b:=\min \{a, b\}$ and $a\vee b:=\max\{a, b\}$. The same
form of the estimates in the case $\alpha=2$ (and $d\ge 3$) were
obtained much earlier in \cite{wid} and \cite{zha}
for the Brownian motion case.

The proofs of \eqref{e:green-for-stable} for symmetric $\alpha$-stable
processes relied heavily on the explicit
formulae for the Green functions and the Poisson kernels of the
ball. Moving away from stable processes, such formulae were not
available and new methods had to be developed. \cite{R}
studied the relativistic $\alpha$-stable process (with relativistic
mass $m>0$) whose characteristic exponent is given by
$\Phi(\theta)=(|\theta|^2+m^{2/\alpha})^{\alpha/2}-m$ and
infinitesimal generator is given by
$m-(-\Delta+m^{2/\alpha})^{\alpha/2}$, and showed that the Green
function of this process in any bounded $C^{1,1}$ domain $D$
satisfies the same sharp estimates \eqref{e:green-for-stable}. Soon
after, \cite{CS3}, using a perturbation method,
established a general result which includes the main result of
\cite{R} as a special case. For different generalizations of the
main result of \cite{R}, see the recent papers \cite{GR, KL}.

Quite recently, \cite{CKS} studied the L\'evy
process which is the sum of independent symmetric $\beta$-stable and
$\alpha$-stable processes, $0<\beta<\alpha<2$. The characteristic
exponent of this process is given by
$\Phi(\theta)=
|\theta|^{\alpha}+|\theta|^{\beta}$ and the
infinitesimal generator by
$-(-\Delta)^{\alpha/2}-(-\Delta)^{\beta/2}$. Sharp two-sided
estimates on the heat kernel of this process in $C^{1, 1}$ open sets
were established in \cite{CKS}. As a by-product of the heat kernel
estimates, sharp Green function estimates of the process in any
bounded $C^{1,1}$ open set were obtained in \cite{CKS}.
Estimates have the form \eqref{e:green-for-stable}. In contrast with
the relativistic stable processes, these Green function estimates
cannot be obtained using the methods of \cite{CS3, GR, KL, R}. The
case $\alpha=2$ (i.e., one of the processes is a Brownian motion)
was covered in \cite{CKSV2} with analogous estimates.

The common feature of these Green function estimates is that both
the distance between the points, $|x-y|$, and distances to $D^c$,
$\delta_D(x),  \delta_D(y)$, appear as arguments of \emph{power
functions}. However, it follows from \cite[Chapter 5]{BBKRSV} that
the asymptotic behavior of the free Green function $G(x,y)$ of many
transient symmetric L\'evy processes is of the form
$$
G(x,y)\sim \frac1{|x-y|^{d-\alpha}\ell(|x-y|^{-2})}\quad \mbox{ as }
|x-y|\to 0
$$
where $\alpha\in (0, 2)$ and $\ell$ is a nontrivial slowly varying
function at infinity. (See, also Theorem \ref{t:Gorigin} below.)
Therefore, Green function estimates of the form
\eqref{e:green-for-stable} cannot be true for these general
symmetric L\'evy processes. The purpose of this paper is to
establish sharp two-sided Green function estimates for these more
general L\'evy processes in open sets of $\R^d$. In our estimates,
$\delta_D(x),  \delta_D(y)$ and $|x-y|$ appear as arguments of
regularly varying functions, not necessarily power functions. In
order to explain our setting and results, let us first note that
stable processes, relativistic stable processes and sums of
independent stable processes can be obtained as subordinate Brownian
motions. Indeed, let $W=(W_t=(W^1_t, \dots, W^d_t):\, t\ge 0)$ be a
$d$-dimensional Brownian motion, and let $S=(S_t:\, t\ge 0)$ be an
independent subordinator. Recall that a subordinator is an
increasing L\'evy process on $[0,\infty)$, which can be
characterized through its Laplace exponent $\phi$:
$\E[\exp\{-\lambda S_t\}]=\exp\{-t\phi(\lambda)\}$, $\lambda>0$. The
process $X=(X_t:\, t\ge 0)$ defined by $X_t:=W_{S_t}$ is called a
subordinate Brownian motion. The infinitesimal generator of $X$ is
$-\phi(-\Delta)$. By choosing the  Laplace exponent $\phi(\lambda)$
as $\lambda^{\alpha/2}$, $(\lambda+m^{2/\alpha})^{\alpha/2}-m$ and
$\lambda^{\alpha/2}+\lambda^{\beta/2}$
respectively, the resulting
subordinate Brownian motion turns out to be a symmetric
$\alpha$-stable process, a relativistic stable process and an
independent sum of $\beta$ and $\alpha$-stable processes
respectively. The Laplace exponent of a subordinator is a Bernstein
function and hence has the representation
$$
\phi(\lambda)=b\lambda +\int_{(0,\infty)}(1-e^{-\lambda t})\, \mu(dt)\, ,
$$
where $b\ge 0$ and $\mu$ is a measure (called the L\'evy measure of
$\phi$) such that $\int_{(0,\infty)}(1\wedge t)\, \mu(dt)<\infty$.
If the measure $\mu$ has a completely monotone density, the Laplace
exponent $\phi$ is called a complete Bernstein function. The
common feature of the Laplace exponents
$\phi(\lambda)=\lambda^{\alpha/2}$,
$\phi(\lambda)=(\lambda+m^{2/\alpha})^{\alpha/2}-m$ and
$\phi(\lambda)=\lambda^{\beta/2}+\lambda^{\alpha/2}$ is that all
three of them are complete Bernstein functions whose behavior at
infinity is given by $\lim_{\lambda\to
\infty}\phi(\lambda)/\lambda^{\alpha/2}=1$. We will see that those
two properties (the latter slightly weakened) are the determining
factors for the Green function estimates \eqref{e:green-for-stable}.

Recall that an open set $D$ in $\bR^d$ ($d\ge 2$) is said to be
a $C^{1,1}$ open set if there exist a localization radius $R>0$ and
a constant $\Lambda>0$ such that for every $z\in\partial D$, there
exist a $C^{1,1}$-function $\psi=\psi_z: \bR^{d-1}\to \bR$
satisfying $\psi (0)=0$,  $\nabla\psi (0)=(0, \dots, 0)$, $\| \nabla
\psi \|_\infty \leq \Lambda$, $| \nabla \psi (x)-\nabla \psi (z)|
\leq \Lambda |x-z|$, and an orthonormal coordinate system $CS_z$:
$y=(y_1, \cdots, y_{d-1}, y_d):=(\wt y, \, y_d)$ with origin at $z$
such that
$$
B(z,R)\cap D=\{ y=(\wt y, \, y_d)\in B(0, R) \mbox{ in } CS_z: y_d
>\psi (\wt y) \}.
$$
The pair $(R, \Lambda)$ is called the characteristics of the
$C^{1,1}$ open set $D$. We remark that in some literature, the
$C^{1,1}$ open set defined above is called a {\it uniform} $C^{1,1}$
open set since $(R, \Lambda)$ is universal for all $z\in \partial D$.
By a
$C^{1, 1}$ open set in $\bR$ we mean an open set which can be
written as the union of disjoint intervals so that the minimum of
the lengths of all these intervals is positive and the minimum of
the distances between these intervals is positive.
Note that a bounded $C^{1,1}$ open set can be disconnected.

The main result of this paper is the following sharp Green function
estimates. In the statement and throughout the paper we use notation
$f(t)\asymp g(t)$ as $t \to \infty$ (resp. $t\to 0+$) if the
quotient $f(t)/g(t)$ stays bounded between two positive constants as
$t\to \infty$ (resp. $t\to 0+$).

\begin{thm}\label{t:Gest2}
Suppose that $X=(X_t:\, t\ge 0)$ is a L\'evy process  whose
characteristic exponent is given by $\Phi(\theta)=\phi(|\theta|^2)$,
$\theta\in \R^d$, where
\begin{equation}\label{e:form-of-phi}
\phi:(0,\infty)\to [0,\infty) \text{  is a complete Bernstein
function such that   } \phi(\lambda)\asymp
\lambda^{\alpha/2}\ell(\lambda)\, , \quad\lambda \to \infty\, ,
\end{equation}
$\alpha\in (0, 2\wedge d)$ and $\ell:(0,\infty)\to (0,\infty)$ is a
measurable, locally bounded function which is slowly varying  at
infinity. When $d\le 2$, we assume an additional assumption, see
\eqref{e:ass4trans}. Then for every bounded $C^{1,1}$ open set $D$
in $\R^d$ with characteristics $(R, \Lambda)$,  there exists
$C_1=C_1(\text{diam}(D), R, \Lambda, \alpha, \ell, d)>1$ such that
the Green function $G_D(x,y)$ of $X$ in $D$ satisfies the following
estimates:
\begin{eqnarray}
\lefteqn{C_1^{-1}\,\left(1 \wedge\frac{(\delta_D(x) \delta_D(y))^{\alpha/2}
\ell(|x-y|^{-2})  }{\sqrt{\ell((\delta_D(x)  )^{-2})\ell((\delta_D(y)
)^{-2})}\, |x-y|^{\alpha}} \right)
\frac1{\ell(|x-y|^{-2})|x-y|^{d-\alpha}}
\nn}\\
&\le & G_D(x,y) \le
C_1\,\,\left(1 \wedge\frac{(\delta_D(x) \delta_D(y))^{\alpha/2}
\ell(|x-y|^{-2})  }{\sqrt{\ell((\delta_D(x)  )^{-2})\ell((\delta_D(y)
)^{-2})}\, |x-y|^{\alpha}} \right)
\frac1{\ell(|x-y|^{-2})|x-y|^{d-\alpha}}.\label{e:Gest21}
\end{eqnarray}
\end{thm}

\medskip

Since the subordinate Brownian motion $X$ is completely determined
by the Laplace exponent $\phi$, one would expect that the above
estimates can be expressed in terms of the function $\phi$ only. And
indeed, an alternative form of \eqref{e:Gest21} reads as follows:
There exists $c>1$ such that
\begin{eqnarray}\label{e:Gest21-alt1}
\lefteqn{c^{-1}\left(1 \wedge
\frac{\phi(|x-y|^{-2})}{\sqrt{\phi(\delta_D(x)^{-2})
\phi(\delta_D(x)^{-2})}}\right)\, \frac{1}{|x-y|^d\,
\phi(|x-y|^{-2})}}\nonumber \\
&\le &G_D(x,y) \le c \left(1 \wedge
\frac{\phi(|x-y|^{-2})}{\sqrt{\phi(\delta_D(x)^{-2})
\phi(\delta_D(x)^{-2})}}\right)\, \frac{1}{|x-y|^d\,
\phi(|x-y|^{-2})}
\end{eqnarray}
(see \eqref{e:Gest21-alt2} below for yet another alternative form of
these estimates). To the best of our knowledge, the above Green
function estimates include all known Green function estimates of
pure-jump transient subordinate Brownian motions in bounded
$C^{1,1}$ open set $D$ in $\R^d$ as special cases. The estimates
above include a lot more processes: In the case $d\ge 3$, our
estimates are valid for all subordinate Brownian motions satisfying
\eqref{e:form-of-phi}. For more concrete examples, see Example
\ref{rexample}.

Let us give the main ingredients of the proof of Theorem
\ref{t:Gest2}. The groundwork has been laid down in the recent paper
\cite{KSV} where a similar class of subordinate Brownian motions was
studied. One difference to the current setting was that in
$\cite{KSV}$ the Laplace exponent was assumed to be precisely
regularly varying at infinity and not just comparable to a
regularly varying function. Another difference is that \cite{KSV}
contains some additional assumptions that we showed to be redundant.
The results of \cite{KSV} are reproved in \cite{KSV2} under
conditions valid in this paper (with the redundant assumptions
removed). When referring to those results we will quote both
sources. The main result of \cite{KSV, KSV2} is the boundary Harnack
principle for nonnegative harmonic functions of the subordinate
Brownian motion $X$ in bounded $\kappa$-fat open sets. Based on the
boundary Harnack principle and using the well-established
methodology of
\cite{B3, H}, we will first obtain Green
function estimates of the form \eqref{e:Gest} (in the spirit of
\cite{B3, H}) in bounded $\kappa$-fat open sets.

Recall  from \cite{SW} that an open set $D$ in $\R^d$ is
$\kappa$-fat if there exists $R_1>0$ such that for each $Q \in
\partial D$ and $r \in (0, R_1)$, $D \cap B(Q,r)$ contains a ball
$B(A_r(Q),\kappa r)$. The pair $(R_1, \kappa)$ is called the
characteristics of the $\kappa$-fat open set $D$. All Lipschitz
domains, non-tangentially accessible domains and John domains are
$\kappa$-fat (cf.~\cite{K2, SW} and the references therein). In
general, the boundary of a $\kappa$-fat open set can be
nonrectifiable.

\begin{thm}\label{t:Gest}
Suppose that $X=(X_t:\, t\ge 0)$ is a L\'evy process satisfying the
same conditions as in Theorem \ref{t:Gest2} and that $D$ is a
bounded $\kappa$-fat open set with characteristics $(R_1, \kappa)$.
Then there exists $C_2=C_2(\text{diam}(D), R_1, \kappa, \alpha,
\ell, d)>1$ such that for every $x, y \in D$,
\begin{equation}\label{e:Gest}
C_2^{-1}\,\frac{g(x) g(y)}{g(A)^2|x-y|^{d-\alpha}\ell(|x-y|^{-2})}\,\le\,
G_D(x,y) \,\le\, C_2\,\frac{g(x) g(y)}{g(A)^2|x-y|^{d-\alpha}\ell(|x-y|^{-2})},
\quad  A \in \BB(x,y),
\end{equation}
where $g$ and $\BB(x,y)$ are defined in \eqref{d:gz0} and
\eqref{d:gz1} respectively.
\end{thm}

In the case $0<c_1 \le \ell (\lambda) \le  c_2<\infty$ for large
$\lambda$, using the Harnack inequality and the boundary Harnack
principle, the above form of Green function estimates has been
established by several authors in special cases. See \cite[Theorem
1.1]{CKSV2}, \cite[Theorem 2.4]{H}  and \cite[Theorem 1]{J}.

To obtain the interior estimates  in Theorem \ref{t:Gest} (i.e.~for
points $x,y$ away from the boundary), we use the asymptotic behavior
of the Green function of $X$ in $\RR^d$  proved in \cite[Theorem
3.2]{KSV2} (see also \cite[Theorem 3.1]{KSV}). Using the interior
estimates and following  the method developed in
\cite{B3, H},
we obtain the full estimates in a bounded $\kappa$-fat open set
using the boundary Harnack principle from \cite{KSV, KSV2}.

Even though a lot of complications occur due to the appearance of
the slowly varying function $\ell$, the proof of Theorem
\ref{t:Gest} is still routine. But the precise estimates
\eqref{e:Gest21} in bounded $C^{1,1}$ open sets are very delicate.
One of the ingredients comes from the fluctuation theory of
one-dimensional L\'evy processes. Let $Z=(Z_t:\, t\ge 0)$ be the
one-dimensional subordinate Brownian motion defined by
$Z_t:=W^d_{S_t}$, and let $V$ be the renewal function of the ladder
height process of $Z$. The function $V$ is harmonic for the process
$Z$ killed upon exiting $(0,\infty)$, and the function $w(x):=
V(x_d){\bf 1}_{\{x_d >0\}}$, $x\in \R^d$, is harmonic for the
process $X$ killed upon exiting the half space $\R^d_+:=\{ x=(x_1,
\dots, x_{d-1}, x_d) \in \bR^d: x_d > 0 \}$ (Theorem \ref{t:Sil}).
Therefore, $w$ gives the correct rate of decay of harmonic
functions near the boundary of $\R^d_+$. This shows the importance
of the fluctuation theory (of one-dimensional L\'evy processes) in
our approach.

The second ingredient is the ``test function'' method applied to the
operator $\mathcal A$ defined by
$$
\mathcal A f(x) =\lim_{\varepsilon \to 0}\int_ {\{y  \in
\R^d:|y-x|>\varepsilon\}}(f(y)-f(x))\, J(y-x)\, dy\, ,
$$
with the domain consisting of functions $f$ for which the limit
exists and is finite. Here $J$ denotes the density of the L\'evy
measure of $X$. On the space of smooth functions with compact
support, this operator coincides with the infinitesimal generator of
$X$. We emphasize  that, compared to the test function methods of
\cite{BBC, CKSV, G}, there are several differences in our approach.
In \cite{BBC, CKSV, G}, appropriate subharmonic and superharmonic
functions of $X$ (or the truncated version of $X$) are chosen as
test functions, first in the case of half spaces and then for $C^{1,
1}$ open sets, and the values of the generator acting on these test
functions are computed in detail. Then suitable combinations of the
test functions are used to find the correct exit distribution
estimates. In \cite{BBC, CKSV, G}, the test functions are power
functions of the form $x \to (x_d)^p$ and the densities of the
L\'evy measures of the processes have closed forms. However, the
 density $J$ of the L\'evy measure of our process does not have a
simple form. We do not even know the  asymptotic behavior of $J$
near infinity in general. Furthermore, in general, power functions
of the form $x \to (x_d)^p$ are neither subharmonic nor
superharmonic functions for our processes, and it is not clear what
are the appropriate choices for the test functions.

Due to the above differences and difficulties, obtaining the correct
boundary decay rate of the Green function in $C^{1,1}$ open set $D$
requires new ideas and approaches. In this paper, we will use the
function $w$ which is smooth and harmonic on the half space, as our
only test function. Using this and the characterization of harmonic
functions recently established in \cite{C}, we show that $\mathcal A
w \equiv 0$ on the half space (Theorem \ref{c:Aw=0}). With this, we
prove the following fact in Lemma \ref{L:Main}, which is the key to
the proof of Theorem \ref{t:Gest2}: If $D$ is a $C^{1,1}$ open set
with characteristics $(R, \Lambda)$, $Q\in \partial D$ and
$h(y)=V(\delta_D(y)) \1_{D\cap B(Q,R)}$, then $\mathcal A h(y)$ is
well defined and bounded for $y\in D$ close enough to the boundary
point $Q$. Using this lemma, we give certain exit distribution
estimates  in Lemma \ref{L:2}, which provide the correct rate of
decay of Green functions near the boundary of $D$. Unlike \cite{BBC,
CKSV, G},  in Lemma \ref{L:2} we do not construct subharmonic and
superharmonic functions on $C^{1,1}$ open set $D$. Instead we use
Dynkin's formula on $h$ to obtain the desired exit distribution
estimates directly. In fact, our approach is  simpler than the
previous approaches and may be used for other types of jump
processes. We hope our approach will shed new light on the
understanding of the boundary behavior of nonnegative harmonic
functions of general Markov processes.

The estimates \eqref{e:Gest21} are best understood in terms of the
renewal function $V$ which provides the exact rate of decay of $G_D$
near the boundary. Let $G$ be the Green function of $X$ in the whole
space $\R^d$. An equivalent form of \eqref{e:Gest21} is given by
\begin{equation}\label{e:Gest21-alt2}
c^{-1}\left(1 \wedge\frac{V(\delta_D(x))V(\delta_D(y))}{V(|x-y|)^2} \right)
G(x, y) \le G_D(x,y)
\le c \left(1 \wedge\frac{V(\delta_D(x))V(\delta_D(y))}{V(|x-y|)^2} \right)
G(x, y).
\end{equation}

By combining the sharp estimates of the Green function in a bounded
$C^{1,1}$ open set with the boundary Harnack principle proved in
\cite{KSV, KSV2} (see Theorem \ref{BHP} below), we  obtain a
boundary Harnack principle with explicit decay rate. In the next
theorem we give an extension to unbounded $C^{1,1}$ open sets.
Recall that, given $Q\in \partial D$, a function $u:\bR^d\to \R$ is
said to vanish continuously on $ D^c \cap B(Q, r)$ if $u=0$ on $ D^c
\cap B(Q, r)$ and $u$ is continuous at every point of $\partial
D\cap B(Q,r)$.

\begin{thm}\label{t:bhp}
Suppose that $X=(X_t:\, t\ge 0)$ is a L\'evy process satisfying the
same conditions as in Theorem \ref{t:Gest2} and that $D$ is a
(possibly unbounded) $C^{1, 1}$ open set in $\bR^d$ with
characteristics $(R, \Lambda)$. Then there exists $C_3=C_3(R,
\Lambda, \alpha, \ell, d)>0$  such that for $r \in (0, (R \wedge
1)/4]$, $Q\in \partial D$ and any nonnegative function $u$ in $\R^d$
that is harmonic in $D \cap B(Q, r)$ with respect to $X$ and
vanishes continuously on $ D^c \cap B(Q, r)$, we have
\begin{equation}\label{e:bhp_m}
\frac{u(x)}{\delta_D(x)^{\alpha/2}\sqrt{\ell((\delta_D(y)
)^{-2})}}\,\le C_3\,
\frac{u(y)}{\delta_D(y)^{\alpha/2}\sqrt{\ell((\delta_D(x)  )^{-2})}}
\qquad \hbox{for every } x, y\in  D \cap B(Q, r/2).
\end{equation}
\end{thm}

An alternative form of \eqref{e:bhp_m} reads as follows: There
exists a constant $c>0$ such that
\begin{equation}\label{e:bhp-m-alt}
\frac{u(x)}{V(\delta_D(x))}\le c \frac{u(y)}{V(\delta_D(y))} \qquad
\hbox{for every } x, y\in  D \cap B(Q, r/2).
\end{equation}
Note that unlike the usual form of the boundary Harnack principle
where one considers the ratio of two harmonic functions, functions
in the denominator of \eqref{e:bhp_m} and \eqref{e:bhp-m-alt} are
not harmonic. Instead, they provide the correct boundary decay of
non-negative harmonic functions. Indeed, an equivalent form of
Theorem \ref{t:bhp} says that there exists a constant $c>1$ such
that for any nonnegative function $u$ in $\R^d$ that is harmonic in
$D \cap B(Q, r)$ with respect to $X$ and vanishes continuously on $
D^c \cap B(Q, r)$ it holds that
$$
c^{-1} V(\delta_D(x)) \le u(x) \le c V(\delta_D(x)) \qquad \hbox{for
every } x\in  D \cap B(Q, r/2).
$$

This paper is organized as follows: In the next section we establish
the setting and notation, prove several new results for complete
Bernstein functions, and describe some of  the known results from
\cite{KSV, KSV2}. In Section 3 we prove the Green function estimates
in bounded $\kappa$-fat open sets. Section 4 is devoted to the Green
function estimates in bounded $C^{1,1}$ open sets.

We will use the following conventions in this paper. The values of
the constants $C_{1},C_2, \dots$, $M$, $\eps_1$ and $R, R_1, R_2,
\cdots$ will remain the same throughout this paper, while $c, c_0,
c_1, c_2, \cdots$ and  $r, r_0, r_1, r_2,\dots $ stand for constants
whose values are unimportant and which may change from one
appearance to another. All constants are positive finite numbers.
The labeling of the constants $c_0, c_1, c_2, \cdots$ starts anew in
the statement of each result. The dependence of the constants  on
dimension $d$ may not be mentioned explicitly. We will use ``$:=$"
to denote a definition, which is  read as ``is defined to be".
Further, $f(t) \sim g(t)$, $t \to 0$ ($f(t) \sim g(t)$, $t \to
\infty$, respectively) means $ \lim_{t \to 0} f(t)/g(t) = 1$
($\lim_{t \to \infty} f(t)/g(t) = 1$, respectively). For any open
set $U$, we denote by $\delta_U (x)$ the distance between $x$ and
the complement of $U$,  i.e., $\delta_U(x)=\text{dist} (x, U^c)$. We
will use $\partial$ to denote the cemetery point and for every
function $f$, we extend its definition to $\partial$ by setting
$f(\partial )=0$. For every function $f$, let $f^+:=f \vee 0$. We
will use $dx$ to denote the Lebesgue measure in $\bR^d$. For a Borel
set $A\subset \bR^d$, we also use $|A|$ to denote its Lebesgue
measure and ${\rm diam}(A)$ to denote the diameter of the set $A$.

\section{Preliminaries}

In this section we collect and explain preliminary results necessary
for further development in Sections 3 and 4. Most of these results
originate from \cite{KSV} where they were proved under somewhat
stronger conditions than in this paper. Their extensions to the
current setting, in particular Theorems \ref{t:Gorigin},
\ref{t:Jorigin} and \ref{BHP}, are given with full proofs in
\cite{KSV2}. Here we prove only results that have not appeared in
\cite{KSV}. Lemma \ref{l:H2-valid} and Propositions \ref{p:chi is
cbf} and \ref{p:chiphi} about complete Bernstein functions may be of
independent interest. The difference between the assumptions in
\cite{KSV} and this paper is discussed in Remark
\ref{r-interpretation-H}.

A $C^{\infty}$ function $\phi:(0,\infty)\to [0,\infty)$ is called a
Bernstein function if $(-1)^n D^n \phi\le 0$ for every positive
integer $n$. Every Bernstein function has a representation $
\phi(\lambda)=a+b\lambda +\int_{(0,\infty)}(1-e^{-\lambda t})\,
\mu(dt) $ where $a,b\ge 0$ and $\mu$ is a measure on $(0,\infty)$
satisfying $\int_{(0,\infty)}(1\wedge t)\, \mu(dt)<\infty$; $a$ is
called the killing coefficient, $b$ the drift and $\mu$ the L\'evy
measure of the Bernstein function. A Bernstein function $\phi$ is
called a complete Bernstein function if the L\'evy measure $\mu$ has
a completely monotone density $\mu(t)$, i.e., $(-1)^n D^n \mu\ge 0$
for every non-negative integer $n$. Here and below, by abuse of
notation  we will denote the L\'evy density by $\mu(t)$. For more on
Bernstein and complete Bernstein functions we refer the readers to
\cite{SSV}.

First, we show that $\phi$ being a complete Bernstein function
implies that its L\'evy density cannot decrease too fast in the
following sense:

\begin{lemma}\label{l:H2-valid}
Suppose that $\phi$ is a complete Bernstein function with L\'evy
density $\mu$. Then there exists $C_4>1$ such that $\mu(t)\le C_4
\mu(t+1)$ for every  $t>1$.
\end{lemma}

\pf Since $\mu$ is a completely monotone function, by Bernstein's
theorem (\cite[Theorem 1.4]{SSV}), there exists a measure $m$ on
$[0,\infty)$ such that $\mu(t)=\int_{[0,\infty)}e^{-tx} m(dx).$
Choose $r>0$ such that $\int_{[0, r]}e^{-x}\, m(dx)\ge \int_{(r,
\infty)}e^{-x}\, m(dx).$ Then, for any $t>1$, we have
\begin{eqnarray*}
\int_{[0, r]}e^{-t x}\, m(dx)&\ge&e^{-(t -1)r} \int_{[0, r]}e^{-x}\, m(dx)\\
&\ge &e^{-(t -1)r}\int_{(r, \infty)}e^{-x}\,  m(dx)\,\ge \,
\int_{(r, \infty)}e^{-t x}\, m(dx).
\end{eqnarray*}
Therefore, for any $t>1$,
$$
\mu(t+1)\ge \int_{[0, r]}e^{-(t+1) x}\, m(dx)\ge  e^{-r}\int_{[0,
r]}e^{- t x}\, m(dx) \ge \frac12\, e^{-r}\int_{[0, \infty)}e^{-t
x}\, m(dx)=\frac12\, e^{-r}\mu(t).
$$
\qed

Suppose that $S=(S_t: t\ge 0)$ is a subordinator with Laplace exponent $\phi$, that is
$$
\E\left[e^{-\lambda S_t}\right]=e^{-t\phi(\lambda)}, \qquad \forall \ t, \lambda>0.
$$
The Laplace exponent of a subordinator is always a Bernstein function.
Let $U(A):=\E\int_0^{\infty}
\1_{\{S_t\in A\}}
\, dt$ denote the potential measure of $S$. If $\phi$ is a complete Bernstein function with infinite L\'evy measure, then the potential measure $U$ has a completely monotone density $u(t)$ (see, e.g., \cite[Remark 10.6 and Corollary 10.7]{SSV}).

Recall that a function $\ell:(0,\infty)\to (0,\infty)$ is slowly varying at infinity if
$$
\lim_{t\to \infty}\frac{\ell(\lambda t)}{\ell(t)}=1\, ,\quad \textrm{for every }\lambda >0\, .
$$

In the remainder of this paper we assume that $\phi$ is a complete Bernstein function and we
will always impose the following

\medskip
\noindent {\bf Assumption (H):} There exist $\alpha\in (0,2)$ and a
function $\ell:(0,\infty)\to (0,\infty)$ which is measurable,
locally bounded and slowly varying at infinity such that
\begin{equation}\label{e:reg-var}
\phi(\lambda) \asymp \lambda^{\alpha/2}\ell(\lambda)\, ,\quad \lambda \to \infty\, .
\end{equation}

\begin{remark}
\label{r-interpretation-H}{\rm
(a) The precise interpretation of \eqref{e:reg-var} will be as follows:
There exists a positive constant $c>1$ such that
$$
c^{-1}\le \frac{\phi(\lambda)}{\lambda^{\alpha/2}\ell(\lambda)}
\le c \qquad \textrm{for all }\lambda \in [1,\infty)\, .
$$
The choice of the interval $[1,\infty)$ is, of course, arbitrary.
Any interval $[a,\infty)$ would do, but with a different constant.
This follows from the assumption that $\ell$ is locally bounded.
Moreover, by choosing $a>0$ large enough, we could dispense with
the local boundedness assumption. Indeed, by \cite[Lemma 1.3.2]{BGT},
every slowly varying function at infinity is locally bounded on
$[a,\infty)$ for $a$ large enough.
Although the choice of interval $[1,\infty)$ is arbitrary, it will
have as a consequence the fact that all relations of the type
$f(t)\asymp g(t)$ as $t \to \infty$ (respectively $t\to 0+$)
following from \eqref{e:reg-var} will be interpreted as
$\tilde{c}^{-1} \le f(t)/g(t) \le \tilde{c}$ for $t\ge 1$
(respectively $0<t\le 1$) for an appropriate constant $\tilde{c}$.

\noindent
(b) The assumption {\bf (H)} is an assumption about the behavior of
$\phi$ at infinity. We make no assumption on $\phi$ near zero.
As a consequence, we will be able to obtain information about
the small scale behavior of the subordinate process,  but almost
nothing can be inferred about its large scale behavior.

\noindent
(c) The main assumption in \cite{KSV} was that $\phi$ is a complete
Bernstein function such that
\begin{equation}\label{e:old-assumption}
\phi(\lambda)=\lambda^{\alpha/2}\ell(\lambda)\, ,\quad \text{for all }
\lambda >0\, ,
\end{equation}
where $\alpha\in (0,2)$ and $\ell$ is a slowly varying function at
infinity. This assumption allows us to obtain exact asymptotic behavior
of various functions. More precisely, some of the results in \cite{KSV}
were of the form $f(t)\sim g(t)$, while with the assumption
\eqref{e:reg-var} we can obtain only the corresponding results in the
weaker form $f(t)\asymp g(t)$. Proofs of these weaker results can be
found in \cite{KSV2}. We note that our current assumptions are indeed
strictly weaker than the ones in \cite{KSV}: There exists a complete
Bernstein function satisfying \eqref{e:reg-var} which is \emph{not}
regularly varying at infinity, see \cite[Example 2.8]{KSV2}.

\noindent
(d) We briefly comment on the other assumptions from \cite{KSV}
which are now removed. The assumption
{\bf A1} in \cite{KSV}
needed for transience
in case $d\le 2$ is replaced by \eqref{e:ass4trans} below.
The assumptions {\bf A2} and {\bf A3}
in \cite{KSV} used
 for the technical lemma
\cite[Lemma 5.32]{BBKRSV} are redundant - see \cite[Lemma 3.1]{KSV2}.
Assumption
{\bf A4} in \cite{KSV}
 is always valid for complete Bernstein function
as proved here in Lemma \ref{l:H2-valid}.  Finally, the assumption (2.5)
in \cite[Proposition 2.2]{KSV} is no longer needed as \cite[Proposition 2.2]{KSV} is now
replaced by Proposition \ref{p:chiphi} below.
}
\end{remark}

It follows from \eqref{e:reg-var} that $\lim_{\lambda\to
\infty}\phi(\lambda)/\lambda=0$ and $\lim_{\lambda\to
\infty}\phi(\lambda)=\infty$, implying that $\phi$ has
no drift and its L\'evy measure is infinite.
Therefore, the potential measure $U$ of the corresponding
subordinator $S$ has a completely monotone density $u$.

The behavior of $u(t)$ and the density $\mu(t)$ of the
L\'evy measure can be inferred from the following result.

\begin{prop}\label{p:zahle}
{\rm(\cite[Theorem 7]{Z})}
Suppose that $\psi$ is a completely monotone function given by
$$
\psi(\lambda)=\int^\infty_0e^{-\lambda t} f(t)\, dt,
$$
where $f$ is a nonnegative decreasing function. Then
$$
f(t)\le \left(1-e^{-1}\right)^{-1} t^{-1}\psi(t^{-1}), \quad t>0.
$$
If, furthermore, there exist $\delta\in (0, 1)$ and $a, t_0>0$ such
that
\begin{equation}\label{e:zahle}
\psi(r \lambda)\le a r^{-\delta} \psi(\lambda), \quad r\ge 1, t\ge 1/t_0,
\end{equation}
then there exists $C_5=C_5(w,f,a,t_0, \delta)>0$ such that
$$
f(t)\ge C_5 t^{-1}\psi(t^{-1}), \quad t\le t_0.
$$
\end{prop}

\medskip

We first apply the above proposition to $\psi(\lambda)=
\phi(\lambda)^{-1}=\int_0^{\infty}e^{-\lambda t}u(t)\, dt$
to obtain the behavior of $u$ near zero:
\begin{equation}\label{e:behofu}
u(t)\asymp t^{-1}\phi(t^{-1})^{-1}\asymp\frac{t^{\alpha/2-1}}
{\ell(t^{-1})}\, , \quad t \to 0+\, .
\end{equation}
Condition \eqref{e:zahle} follows from \eqref{e:reg-var} by
use of Potter's theorem (cf. \cite[Theorem 1.5.6]{BGT}).
By applying \eqref{e:behofu} to the complete Bernstein function
$\lambda\mapsto \lambda/\phi(\lambda)$ (\cite[Proposition 7.1]{SSV})
one obtains the following behavior of $\mu(t)$ near zero:
\begin{equation}\label{e:behofmu}
\mu(t)\asymp t^{-1}\phi(t^{-1})\asymp t^{-\alpha/2-1}
\ell(t^{-1})\, , \quad t\to 0+\, .
\end{equation}
We refer the reader to \cite[Theorem 2.9, Theorem 2.10]{KSV2}
for the detailed proofs of \eqref{e:behofu} and \eqref{e:behofmu}.
The corresponding precise asymptotics are given in \cite[p.~1603]{KSV}
under the assumption \eqref{e:old-assumption}.

A consequence of the asymptotic behavior \eqref{e:behofmu} of $\mu(t)$
is that for any $K>0$ there exists $c=c(K) >1$ such that
\begin{equation}\label{e:mu at zero}
\mu(t)\le c\, \mu(2t), \qquad \, \forall t\in (0, K).
\end{equation}
The behavior of $\mu(t)$ at infinity has already been determined in
Lemma \ref{l:H2-valid}: There exists a constant $c >1$ such that
\begin{equation}\label{e:mu at infty}
\mu(t)\le c\, \mu(t+1), \qquad \, \forall\,  t>1.
\end{equation}
This property of $\mu$ was assumed in \cite{KSV} as {\bf A4}, but we
have shown in Lemma \ref{l:H2-valid} that it always holds true.

\medskip

We consider now one-dimensional subordinate Brownian motions.
Let $B=(B_t:\, t\ge 0)$ be a Brownian motion in $\R$, independent of $S$, with
$$
\E\left[e^{i\theta(B_t-B_0)}\right]
=e^{-t\theta^2}, \qquad\,  \forall \theta\in \R, t>0.
$$
The subordinate Brownian motion $Z=(Z_t:t\ge 0)$ in $\R$ defined by
$Z_t=B_{S_t}$ is a symmetric L\'evy process with the characteristic exponent
$\Phi(\theta)=\phi(\theta^2)$ for all $\theta\in \R.$

Let $\overline{Z}_t:=\sup\{0\vee Z_s:0\le s\le t\}$ be the supremum process of $Z$ and
let $L=(L_t:\, t\ge 0)$ be a local time of $\overline{Z}-Z$ at $0$. $L$ is
also called a local time of the process $Z$ reflected at the supremum.
The right continuous inverse $L^{-1}_t$ of $L$ is a subordinator and
is called the ladder time process of $Z$. The process $\overline{Z}_{L^{-1}_t}$ is also
a subordinator and is called the ladder height process of $Z$.
(For the basic properties of the ladder time and ladder height
processes, we refer our readers to \cite[Chapter 6]{Ber}.)

Let $\chi$ be the Laplace exponent of the ladder height process of
$Z$. It follows from \cite[Corollary 9.7]{Fris} that
\begin{equation}\label{e:formula4leoflh}
\chi(\lambda)=
\exp\left(\frac1\pi\int^{\infty}_0\frac{\log(\Phi(\lambda\theta))}
{1+\theta^2}d\theta
\right)
=\exp\left(\frac1\pi\int^{\infty}_0\frac{\log(
 \phi(\lambda^2\theta^2))}{1+\theta^2}d\theta
\right), \quad \forall \lambda>0.
\end{equation}
Using \eqref{e:formula4leoflh}, it was proved in \cite[Proposition 2.1]{KSV} that,
if $\phi$ is a special Bernstein
function, so is $\chi$, i.e., $\lambda\mapsto \lambda /\chi(\lambda)$ is also a
Bernstein function. The next result tells us that such relation is also
true for complete Bernstein functions.
For the proof  we will need the following fact, see \cite[Theorem 6.10]{SSV}:
If $\phi$ is a complete Bernstein function,
then there exist a real number $\gamma$
and a $[0, 1]$-valued function $\eta$ on $(0, \infty)$ such that
\begin{equation}\label{e:exp-repr}
\log \phi(\lambda)=\gamma
+\int^\infty_0\left(\frac{t}{1+t^2}-\frac1{\lambda+t}\right)\eta(t)dt.
\end{equation}

\begin{prop}\label{p:chi is cbf}
Suppose $\phi$, the Laplace exponent of the subordinator $S$, is a
complete Bernstein function.
Then the Laplace exponent $\chi$ of the ladder height process of the
subordinate Brownian motion $Z_t=B_{S_t}$
is also  a complete Bernstein function.
\end{prop}

\pf
By  \eqref{e:formula4leoflh} and \eqref{e:exp-repr},  we have
\begin{eqnarray*}
\log \chi(\lambda)
=\frac{\gamma}2+\frac1{\pi}\int^\infty_0\int^\infty_0
\left(\frac{t}{1+t^2}-\frac1{\lambda^2\theta^2+t}\right)\eta(t)dt
\frac{d\theta}{1+\theta^2}.
\end{eqnarray*}
By using $0\le \eta(t)\le 1$, we have
\begin{eqnarray*}
\eta(t)\left|\frac{t}{1+t^2}-\frac1{\lambda^2\theta^2+t}\right|
\frac{1}{1+\theta^2}&\le& \frac{1}{1+t^2} \frac{1}{1+\theta^2}\left(\frac1{\lambda^2
\theta^2+t }+ \frac{\lambda^2 \theta^2t}{\lambda^2 \theta^2+t
}\right)\\
&  \le& \frac{1}{1+t^2} \left(\frac1{\lambda^2
\theta^2+t }+ \frac{\lambda^2 t}{\lambda^2 \theta^2+t
}\right).
\end{eqnarray*}
Since
\begin{eqnarray*}
\int^\infty_0\frac1{\lambda^2 \theta^2+t }\, d\theta
=\frac1t\int^\infty_0\frac1{\frac{\lambda^2\theta^2}t+1}\, d\theta=\frac1t\frac{\sqrt{t}}{\lambda}\int^\infty_0\frac1{\gamma^2+1}\, d\gamma
= \frac{\pi}{2\lambda\sqrt{t}},
\end{eqnarray*}
we can use  Fubini's theorem to get
\begin{eqnarray}
\log \chi(\lambda)
&=&\frac{\gamma}2+\int^\infty_0\left(\frac{t}{2(1+t^2)}-
\frac1{2\sqrt{t}(\lambda+\sqrt{t})}\right)\eta(t)dt \label{e:logexp4chi} \\
&=&\frac{\gamma}2+\int^\infty_0\left(\frac{t}{2(1+t^2)}-
\frac1{2(1+t)}\right)\eta(t)dt+\int^\infty_0
\left(\frac1{2(1+t)}-\frac1{2\sqrt{t}(\lambda+\sqrt{t})}\right)\eta(t)dt\nonumber\\
&=&\gamma_1+\int^\infty_0\left(\frac{s}{1+s^2}-\frac1{\lambda+s}
\right)\eta(s^2)ds\, ,\nonumber
\end{eqnarray}
Applying \cite[Theorem 6.10]{SSV}
we get that $\chi$ is a
complete Bernstein function. \qed
\begin{remark}\label{r:kwasnicki}{\rm
The above result has been independently proved in \cite[Lemma 4]{Kwa}.
}
\end{remark}

The next result relates the behavior of $\chi$ with that of $\phi$.
It will be used to obtain the asymptotic behavior of $\chi$ at infinity.

\begin{prop}\label{p:chiphi}
Suppose that $\phi$, the Laplace exponent of the subordinator $S$,
is a complete Bernstein function. Then the Laplace exponent $\chi$ of the
ladder height process of $Z$ satisfies
\begin{equation}\label{e:chi-and-phi}
e^{-\pi/2}\sqrt{\phi(\lambda^2)}\le \chi(\lambda)\le e^{\pi/2}\sqrt{\phi(\lambda^2)}\, ,\qquad
\textrm{for all }\lambda>0\, .
\end{equation}
\end{prop}

\pf
By the representations \eqref{e:exp-repr}  and \eqref{e:logexp4chi},  we get that for all $\lambda >0$
\begin{eqnarray*}
\left|\, \log \chi(\lambda)-\frac12 \log \phi(\lambda^2)\right| &=&
\frac12\left|\int_0^{\infty}\left(\Big(\frac{t}{1+t^2}- \frac1{\sqrt{t}
(\lambda+\sqrt{t})}\Big) -\Big(\frac{t}{1+t^2}-\frac1{\lambda^2+t}
\Big)\right)\eta(t)\, dt\right|\\
&\le &\frac12\int_0^{\infty}\frac{\lambda(\sqrt{t}+\lambda)}{(\lambda^2+t)
\sqrt{t}(\lambda+\sqrt{t})}\, dt=\frac12\int_0^{\infty}\frac{\lambda}{(\lambda^2+t)\sqrt{t}}\, dt =\frac{\pi}{2}\, .
\end{eqnarray*}
This implies that
$
-\pi/2 \le \log \chi(\lambda)-\frac12 \log \phi(\lambda^2) \le \pi/2$
for every $\lambda>0$, which is \eqref{e:chi-and-phi}
\qed

\begin{remark}\label{r:only-cbf}{\rm
We note that for the last two propositions we only need to assume
that $\phi$ is a complete Bernstein function; the assumption {\bf (H)} is not used.
}
\end{remark}

Let $V$ denote the potential measure of the ladder height process
of $Z$. We will also use $V$ to denote the corresponding renewal
function, $V(t):=V((0,t))$. It follows from \eqref{e:reg-var} and
\eqref{e:chi-and-phi} that $\lim_{\lambda\to \infty}\chi(\lambda)/
\lambda=0$ and $\lim_{\lambda\to \infty}\chi(\lambda)=\infty$.
Therefore, the ladder height process of $Z$ has no drift and
has infinite L\'evy measure. This suffices to conclude that the
potential measure $V$ has a density denoted by $v$, and the
renewal function can be written as $V(t)=\int_0^t v(s)\, ds$.
Since $\chi$ is a complete Bernstein function, $v$ is
completely monotone. We record these facts as

\begin{corollary}\label{c:vsm}
Suppose $\phi$, the Laplace exponent of the subordinator $S$, is a
complete Bernstein function satisfying the assumption {\bf (H)}.
Then the potential measure  of the ladder height process of the
subordinate Brownian motion $Z_t=B_{S_t}$ has
a completely monotone density $v$. In particular, $v$ and the
renewal function $V$
are  $C^{\infty}$ functions.
\end{corollary}

The smoothness of the renewal function $V$ of the ladder height
process  $Z$ will be used later in this paper.

Similarly to the case of the density $u$ of the potential measure $U$ of
the subordinator $S$ in \eqref{e:behofu}, by using Proposition \ref{p:chiphi},
we can obtain the asymptotic behavior of the renewal function $V$ and
its density $v$ of the ladder height process of $Z$:
\begin{eqnarray}
V(t) &\asymp &
\phi(t^{-2})^{-1/2}\asymp\,\frac{t^{\alpha/2}}{(\ell(t^{-2}))^{1/2}}\, ,
\quad t \to 0+\, ,\label{e:behofV}\\
v(t) &\asymp &
t^{-1}\phi(t^{-2})^{-1/2}\asymp\,\frac{t^{\alpha/2-1}}{
(\ell(t^{-2}))^{1/2}}\, ,  \quad t \to 0+\, ,\label{e:behofv}
\end{eqnarray}
see \cite[Proposition 3.9]{KSV2}. The corresponding precise
asymptotics are given in \cite[Proposition 2.7]{KSV} under the
assumption \eqref{e:old-assumption}.

\medskip
We next consider multidimensional subordinate Brownian motions.
Let $W=(W_t=(W^1_{t}, \dots,
W^d_{t}):t\ge 0)$ be a Brownian motion in $\R^d$ with
$$
\E\left[e^{i\theta\cdot(W_t-W_0)}\right]=
e^{-t|\theta|^2}, \qquad \forall\,  \theta \in \R^d, t>0\, ,
$$
and let $S$ be a subordinator independent of $W$ with Laplace exponent $\phi$.
In the remainder of this paper we will use $X=(X_t:\,
t\ge 0)$ to denote the subordinate Brownian motion defined by
$X_t=W_{S_t}$. The process $X$ is a (rotationally) symmetric L\'evy
process with the characteristic exponent given by
$\Phi(\theta)=\phi(|\theta|^2)$.
It is easy to check that when $d\ge 3$ the process $X$ is transient.
This follows from the criterion of Chung-Fuchs type (e.g., \cite[p.~252]{Sa})
which for the subordinate Brownian motion $X$ translates to the following:
$X$ is transient if and only if
$$
\int_{0+} \frac{\lambda^{d/2-1}}{\phi(\lambda)}\, d\lambda < +\infty\, .
$$
Since transience is a global property of the process, it cannot be
inferred from the behavior of $\phi$ at infinity. For example,
$\phi(\lambda)=\log(1+\lambda)+\lambda^{\alpha/2}$, $\alpha\in (0,2)$,
is a complete Bernstein function satisfying \eqref{e:reg-var}, but
the corresponding  subordinate Brownian motion is recurrent in
dimensions 1 and 2. To ensure transience, we will assume that in
the case $d\le 2$, there exists $\gamma\in [0, d/2)$ such that
\begin{equation}\label{e:ass4trans}
\liminf_{\lambda \to 0}\frac{\phi(\lambda)}{\lambda^{\gamma}}>0\, .
\end{equation}
An immediate consequence of this assumption and \cite[Corollary 2.6]{KSV2}
is that the potential density $u$ of $S$ satisfies $u(t)\le c t^{\gamma -1}$
for all $t\ge 1$, where $c>0$ is some positive constant
(cf.~assumption {\bf A1} from \cite{KSV}).

Transience of the process $X$ ensures that the Green function
$G(x,y)$, $x,y\in \R^d$, is well defined. By spatial homogeneity we
may write $G(x,y)=G(x-y)$, where the function
$G$ is radial and given by the
following formula,
\begin{equation}\label{e:representation-G}
G(x)=\int^{\infty}_0(4\pi t)^{-d/2}e^{-|x|^2/(4t)}u(t)dt, \qquad x\in \R^d.
\end{equation}
Since $u$ is decreasing,  we see that $G$ is
radially decreasing and continuous in $\R^d\setminus \{0\}$.

The L\'evy measure of the process $X$ has a density $J$, called
the L\'evy density, given by
$$
J(x)=\int^{\infty}_0(4\pi t)^{-d/2}e^{-|x|^2/(4t)}\mu(t)\, dt, \qquad x\in \R^d.
$$
Thus $J(x)=j(|x|)$ with
\begin{equation}\label{e:representation-j}
j(r)
:=\int^{\infty}_0(4\pi t)^{-d/2}e^{-r^2/(4t)}\mu(t)\, dt, \qquad r>0.
\end{equation}
Note that the function $r\mapsto j(r)$ is continuous and decreasing
on $(0, \infty)$. We will sometimes use the notation $J(x,y)$ for
$J(x-y)$.

We discuss now the behavior of $G$ and $j$ near the origin.
Under the assumption \eqref{e:old-assumption}, the precise
asymptotic behavior was obtained in \cite[Theorem 3.1, Theorem 3.2]{KSV}
by using precise asymptotic behavior of the potential density $u$,
and, respectively, L\'evy density $\mu$. These two results were
proved by use of \cite[Lemma 5.32]{BBKRSV}, which required
additional assumptions which were stated as {\bf A2} and {\bf A3} in
\cite{KSV}. It turned out that by using Potter's theorem
(\cite[Theorem 1.5.6]{BGT}) one can circumvent those assumption
and still obtain the conclusion of the lemma. The details are
provided in \cite[Lemma 3.1]{KSV2}. The lemma combined with
\eqref{e:behofu} (resp.~\eqref{e:behofmu}) and representation
\eqref{e:representation-G} (resp.~\eqref{e:representation-j})
gives the following asymptotic behavior of $G$ (resp.~$J$ under
the assumption \eqref{e:reg-var}:

\begin{thm}\label{t:Gorigin}{\rm(\cite[Theorem 3.2]{KSV2})}
Suppose that the Laplace exponent $\phi$ is a complete Bernstein
function satisfying the assumption
{\bf (H)} and that $\alpha\in (0,2\wedge d)$. In the case $d\le 2$,
we further assume \eqref{e:ass4trans}. Then
$$
G(x)\asymp \frac1{|x|^{d}\phi(|x|^{-2})}\asymp
\frac1{|x|^{d-\alpha}\ell(|x|^{-2})}, \qquad |x|\to 0.
$$
\end{thm}

\begin{remark}\label{cond&zahle}{\rm  Since $\alpha$ is always assumed to be in $(0, 2)$,
the assumption $\alpha\in (0, 2\wedge d)$ in the theorem above makes a difference only
in the case $d=1$.
}
\end{remark}

\begin{thm}\label{t:Jorigin}{\rm(\cite[Theorem 3.4]{KSV2})}
Suppose that the Laplace exponent $\phi$ is a complete Bernstein
function satisfying the assumption {\bf (H)}.
Then
$$
J(x)=j(|x|)\asymp \frac{\phi(|x|^{-2})}{|x|^{d}}\asymp
\frac{\ell(|x|^{-2})}{|x|^{d+\alpha}}, \qquad |x|\to 0.
$$
\end{thm}

Using \eqref{e:mu at zero} and \eqref{e:mu at infty},  and repeating
the proof of \cite[Lemma 4.2]{RSV} we get that

\begin{description}
\item{(a)} For any $K>0$, there exists $c=c(K)
>1$ such that
\begin{equation}\label{H:1}
j(r)\le c\, j(2r), \qquad \forall r\in (0, K).
\end{equation}
\item{(b)} There exists $c
>1$ such that
\begin{equation}\label{H:2}
j(r)\le c\, j(r+1), \qquad \forall r>1.
\end{equation}
\end{description}
Therefore by \cite[Theorem 2.2 and Section 3.1]{SV04} (see also
\cite[Theorem 5.66]{BBKRSV}, \cite[Theorem 4.7, Corollary 4.8]{KSV2} and \cite{RSV})
the Harnack inequality
is valid for the process $X$. Before we state the Harnack inequality,
we recall the definition of harmonic functions.

For any open set $D$, we use $\tau_D$ to denote the first exit
time of $D$, i.e., $\tau_D=\inf\{t>0: \, X_t\notin D\}$.

\begin{defn}\label{def:har1}
Let $D$ be an open subset of $\R^d$.
A function $u$ defined on $\R^d$ is said to be

\begin{description}
\item{\rm{(1)}}  harmonic in $D$ with respect to $X$ if
$\E_x\left[|u(X_{\tau_{B}})|\right] <\infty$ and $
u(x)= \E_x\left[u(X_{\tau_{B}})\right]$,
$ x\in B,
$
for every open set $B$ whose closure is a compact
subset of $D$;

\item{\rm{(2)}}
regular harmonic in $D$ with respect to $X$ if
it is harmonic in $D$ with respect to $X$ and
for each $x \in D$,
$u(x)= \E_x\left[u(X_{\tau_{D}})\right].$
\end{description}
\end{defn}

\begin{thm}[Harnack inequality]\label{T:Har}
There exists $C_{6}>0$ such that, for any $r\in (0, 1/4)$,
$x_0\in\R^d$, and any function $u$ which is nonnegative on $\R^d$
and harmonic with respect to $X$ in $B(x_0, 16r)$, we have
$$
\sup_{y \in B(x_0, r/2)}u(y) \le C_6 \inf _{y \in B(x_0, r/2)}u(y).
$$
\end{thm}

>From now we will always assume that $\phi$ is a complete Bernstein
function satisfying the assumption ({\bf H}) for $\alpha\in
(0, 2\wedge d)$ and the additional \eqref{e:ass4trans} in the case $d\le 2$.
We will no longer explicitly mention these assumptions

For any open set $D$ in $\R^d$, we will use $G_D(x,y)$
to denote the Green function of $X$ in $D$.
Using the continuity
and the radial decreasing property of $G$, we can easily check that
$G_D$ is continuous in $(D\times D)\setminus\{(x, x): x\in D\}$.
We will frequently use the well-known fact that $G_D(\cdot, y)$ is
harmonic in $D\setminus\{y\}$, and regular harmonic in $D\setminus
\overline{B(y,\varepsilon)}$ for every $\varepsilon >0$.

Using the L\'{e}vy system for $X$, we know that for every
bounded open subset $D$ and  every $f \ge 0$ and $x \in D$,
\begin{equation}\label{newls}
\E_x\left[f(X_{\tau_D});\,X_{\tau_D-} \not= X_{\tau_D}  \right]
=  \int_{\overline{D}^c} \int_{D}
G_D(x,z) J(z-y) dz f(y)dy.
\end{equation}

Now we prove the following version of the Harnack inequality for
$X$.

\begin{thm}\label{HP2}
Let $L>0$. There exists a positive constant $C_{7}=C_{7}(L)>0$ such that the
following is true: If $x_{1}, x_{2}\in \R^d$ and $r\in (0,1)$ are such
that $|x_{1}-x_{2}|< Lr$, then for every nonnegative function $u$ which
is harmonic with respect to $X$ in $B(x_{1}, r)\cup B(x_{2},r)$, we
have
$$
C_{7}^{-1}u(x_{2})\,\leq \,u(x_{1}) \,\leq\, C_{7} u(x_{2}).
$$
\end{thm}

\pf
By
\cite[Proposition 4.10]{KSV2}
 (see also \cite[Proposition 3.8]{KSV}),
for every $r\in (0,1)$, every $x\in \R^d$
and  every $y \in B(x, \frac{r}8)$ it holds that
\begin{equation}\label{ddd} K_{B(x,
\frac{r}8)}(x,y):=\int_{B(x, \frac{r}8)} G_{B(x, \frac{r}8)}(x,z)
J(z-y) dz \,\ge\, c_1\, j(|x-y|)\frac{r^\alpha}{\ell(r^{-2})}\, ,
\end{equation}
with some positive constant $c_1>0$.

Now let $r \in (0,1)$,  $x_{1}, x_{2}\in \R^d$ be such that $|x_1-x_2|< Lr$
and let $u$ be a nonnegative function which is harmonic in
$B(x_{1}, r)\cup B(x_{2},r)$ with respect to $X$.

If $ |x_{1}- x_{2}|  <  \frac14 r$, then since $r<1$, the
theorem is true from Theorem \ref{T:Har}. Thus we only need to
consider the case when $\frac14 r \le |x_{1}- x_{2}| \le L r$ with
$L > \frac14$.

Let $w\in B(x_1, \frac{r}{8})$. Because $|x_2-w| \le
|x_1-x_2|+|w-x_1| < (L+\frac18)r \le 2Lr$, by the monotonicity of
$j$ and \eqref{ddd}
\begin{equation}\label{e:asdads1}
K_{B(x_2, \frac{r}8)}(x_2,w)  \,\ge\, c_1\,
j(2Lr)\frac{r^\alpha}{\ell(r^{-2})}\, .
\end{equation}
For any $y\in B(x_1, \frac{r}8)$, $u$ is regular harmonic in $B(y,
\frac{7r}8)\cup B(x_1, \frac{7r}8)$. Since $|y-x_1|< \frac{r}8$, by
Theorem \ref{T:Har}
\begin{equation}\label{e:asdads2} u(y)\ge c_2
u(x_1), \quad y\in B(x_1, \frac{r}8),
\end{equation}
for some constant $c_2>0$. Therefore, by \eqref{newls} and
\eqref{e:asdads1}--\eqref{e:asdads2},
\begin{eqnarray*}
u(x_2) &=& \E_{x_2}\left[u(X_{\tau_{B(x_2,\frac{r}8) }})\right]
  \ge \E_{x_2}\left[u(X_{\tau_{B(x_2,\frac{r}8)}});
X_{\tau_{B(x_2,\frac{r}8)}} \in B(x_1,\frac{r}8) \right]\\
  &\ge &c_2 \,u(x_1)\, \P_{x_2}\left(X_{\tau_{B(x_2,\frac{r}8)}}
\in B(x_1,\frac{r}8) \right)
 = c_2\, u(x_1) \int_{B(x_1,\frac{r}8)} K_{B(x_2,\frac{r}8)} (x_2,w)\, dw \\
 &\ge  &c_3\, u(x_1) \left|B(x_1,\frac{r}8)\right|
  j(2Lr)\frac{r^\alpha}{\ell(r^{-2})}\,= \, c_4\, u(x_1)
 j(2Lr)\frac{r^{\alpha+d}}{\ell(r^{-2})}.
\end{eqnarray*}
Thus, by Theorem \ref{t:Jorigin} and the inequality above, there
exists a constant $c_5>0$ such that for all $r\in (0,1)$,
$u(x_2)  \ge   c_5 u(x_1) \frac{j(2Lr)}{j(r)}$. The
right-hand side is by \eqref{H:1} greater than $c_5 c_6^{\log 2L/ \log 2} u(x_1)$
where $c_6=C_4^{-1}(4L)\in (0,1)$. We have thus proved
the right-hand side inequality in the conclusion of the theorem. The
inequality on the left-hand side can be proved similarly. \qed

\medskip

In \cite{KSV} we have established the boundary Harnack principle under the
assumption \eqref{e:old-assumption} (and the additional transience
assumption in case $d\le2$)
for $\kappa$-fat open sets. Even though we only explicitly stated
the results for $d\ge 2$ in \cite{KSV},
the results and arguments there are in fact valid for $d=1$ also.
Under the current assumptions, the same result is reproved in \cite{KSV2}.

\begin{thm}\label{BHP}\rm{(\cite[Theorem 4.8]{KSV},
\cite[Theorem 4.22]{KSV2})}
Suppose that $D$ is  a bounded  $\kappa$-fat open set with the
characteristics $(R_1, \kappa)$. There exists a  constant
$C_{8}=C_8(d,\alpha,\ell, R_1, \kappa)>1$  such
that if $r \le R_1\wedge \frac14$ and $Q\in\partial D$,
then for any nonnegative functions $u, v$ in $\R^d$
which are regular harmonic in $D\cap B(Q, 2r)$ with respect to $X$
and vanish in $D^c \cap B(Q, 2r)$, we have
$$
C_{8}^{-1}\,\frac{u(A_r(Q))}{v(A_r(Q))}\,\le\, \frac{u(x)}{v(x)}\,\le
C_{8}\,\frac{u(A_r(Q))}{v(A_r(Q))}, \qquad x\in D\cap B(Q, \frac{r}2)\, .
$$
\end{thm}

Before concluding this section, we give some examples satisfying our
assumptions.

\begin{example}\label{rexample}
{\rm Suppose that $0<\beta<\alpha<2$, $0<\gamma < 2-\alpha$ and define
$$
\phi_1(\lambda)=\lambda^{\alpha/2}, \qquad
\phi_2(\lambda)=(\lambda+1)^{\alpha/2}-1, \qquad
\phi_3(\lambda)=\lambda^{\alpha/2}+\lambda^{\beta/2},
$$
$$
\phi_4(\lambda)=\lambda^{\alpha/2}(\log(1+\lambda))^{\gamma/2}\qquad \text{and} \qquad
\phi_5(\lambda)=\lambda^{\alpha/2}(\log(1+\lambda))^{-\beta/2}.
$$
Then $\phi_i$, $i=1,\dots,5$,  are  complete Bernstein functions which can be written
as
$$
\phi_i(\lambda)=\lambda^{\alpha/2}\ell_i(\lambda), \qquad i=1,\dots,5\, ,
$$
with
$$
 \ell_1(\lambda)=1, \qquad
\ell_2(\lambda)=\left((\lambda+1)^{\alpha/2}-1\right)\lambda^{-\alpha/2},
\qquad \ell_3(\lambda)=1+\lambda^{(\beta-\alpha)/2},
$$
$$
\ell_4(\lambda)=(\log(1+\lambda))^{\gamma/2}
\qquad \text{and} \qquad
\ell_5(\lambda)=(\log(1+\lambda))^{-\beta/2}.
$$
As already mentioned in the introduction, the subordinate Brownian
motion corresponding to $\phi_1$ is a
symmetric $\alpha$-stable process, the subordinate Brownian motion
corresponding to $\phi_2$ is a relativistic $\alpha$-stable process
and the subordinate Brownian motion corresponding to $\phi_3$ is the
sum of a symmetric $\alpha$-stable process and an independent
symmetric $\beta$-stable process. The subordinate Brownian motions
corresponding to $\phi_4$ and $\phi_5$ were discussed in
\cite{BBKRSV}.

In the case $d\ge 3$, the only condition on the complete Bernstein
function $\phi$ is \eqref{e:form-of-phi}, so
we can use Proposition 7.1, Corollary 7.9, Propositions 7.10--7.11,
Corollary 7.12 of \cite{SSV}
to come up with infinitely many examples of such functions, e.g.:
(i) $\lambda^{\alpha/2}(\log(1+\log(1+\lambda^{\gamma/2})^{\delta/2}))^{\beta/2}$,
$\alpha, \gamma, \delta \in (0, 2), \beta\in (0, 2-\alpha]$; (ii) $\lambda^{\alpha/2}(\log(1+\log(1+\lambda^{\gamma/2})^{\delta/2}))^{-\beta/2}$,
$\alpha, \gamma, \delta \in (0, 2), \beta\in (0, \alpha]$.

All of the listed example satisfy the stronger condition
\eqref{e:old-assumption}. As already mentioned, \cite[Example 2.8]{KSV2}
provides an example of a complete Bernstein function which satisfies
\eqref{e:reg-var}, but not \eqref{e:old-assumption}.
}
\end{example}

\section{Green function estimates on bounded $\kappa$-fat open sets}\label{s:3}

In this section we will establish sharp two-sided Green function
estimates for $X$
in any bounded $\kappa$-fat open subset $D$ of $\R^d$,
$d\ge 1$.
Our standing assumption is that $\phi$ is a complete Bernstein function
satisfying the assumption {\bf (H)} for $\alpha\in (0, 2\wedge d)$ and the additional
assumption  \eqref{e:ass4trans}
when $d\le 2$.

\begin{lemma}\label{G_G}
There exist $R_2=R_2(\ell) >0$, $L_1>2$  and $C_{9}>0$ such that
$$
G(x) - G(L_1 x) \ge C_{9}\frac1{|x|^{d-\alpha}\ell(|x|^{-2})}
\qquad \text{for every } |x|<R_2.
$$
\end{lemma}
\pf By Theorem \ref{t:Gorigin} there exists a constant $c>1$ such
that for all $x\in \R^d\setminus\{0\}$ with $|x|<1$ it holds that
$$
\frac{c^{-1}}{|x|^{d-\alpha}\ell(|x|^{-2})}\le G(x) \le
\frac{c}{|x|^{d-\alpha}\ell(|x|^{-2})}\, .
$$
Choose $L_1=(4c^2)^{\frac{1}{d-\alpha}}\vee 2$ so that
$c^2/L_1^{d-\alpha}\le 1/4$. Since $\ell$ is slowly varying
at infinity, there exists $r_1<1$ such that
$$
\frac{\ell(|x|^{-2})}{\ell(|L_1 x|^{-2})}\le 2
$$
whenever $0<|x|<r_1$. Let $R_2=r_1 \wedge L_1^{-1}$. Then for
$x\in \R^d\setminus \{0\}$ we have
\begin{eqnarray*}
G(x) - G(L_1 x) & \ge & \frac{c^{-1}}{|x|^{d-\alpha}
\ell(|x|^{-2})}- \frac{c}{|L_1 x|^{d-\alpha}\ell(|L_1 x|^{-2})}\\
&=&\frac{c^{-1}}{|x|^{d-\alpha}\ell(|x|^{-2})}
\left(1-\frac{c^2}{L_1^{d-\alpha}}\,
\frac{\ell(|x|^{-2})}{\ell(|L_1 x|^{-2})}\right)\\
&\ge & \frac{1}{2c} \, \frac{1}{|x|^{d-\alpha}\ell(|x|^{-2})}\, .
\end{eqnarray*}
\qed

\begin{prop}\label{GI}
For every bounded open set $D$, there exists a constant
$C_{10}>1$ such that
 \begin{equation}\label{ub}
G_D(x,y) \le C_{10} \frac1{|x-y|^{d-\alpha}\ell(|x-y|^{-2})}\, ,\qquad
\mbox{ for }x,y \in D\, ,
 \end{equation}
and
\begin{equation}\label{lb}
G_D(x,y) \,\ge\, C_{10}^{-1}\, \frac1{|x-y|^{d-\alpha}\ell(|x-y|^{-2})},
\qquad \mbox{ for }x, y \in D \mbox{ with } L_1|x-y| \le
\delta_D(x) \wedge \delta_D(y),
\end{equation}
where $L_1$ is the constant in Lemma \ref{G_G}.
\end{prop}

\pf
Since $G_D(x,y) \le G(x,y)$, $D$ is bounded and $\ell$ locally bounded,
\eqref{ub} is an immediate consequence of Theorem \ref{t:Gorigin}.
Now we show \eqref{lb}. Without loss of generality, we assume
$\delta_D(y) \le \delta_D(x)$, and let $M:=\mathrm{diam}(D)$.
We consider three cases separately:

\noindent (a) $\delta_D(y) \le R_2$: Since $|x-y| \le
\delta_D(y)/L_1$, $ |X_{\tau_{B(y, \delta_D(y))}}-y| \ge \delta_D(y)
\ge L_1 |x-y| . $ Thus by the monotonicity of $G$ and Lemma
\ref{G_G},
\begin{eqnarray*}
\lefteqn{G_D(x,y) \,\ge\, G_{B(y, \delta_D(y))}(x,y)
\,
=\, G(x,y) - \E_x\left[     G( X_{\tau_{B(y, \delta_D(y))}}, y) \right]}\\
&\ge & G(x-y) - G(L_1(x-y)) \,\ge\, c_1
\frac1{|x-y|^{d-\alpha}\ell(|x-y|^{-2})},  \qquad \forall |x-y| \le
\frac{\delta_D(y)}{L_1}\, .
\end{eqnarray*}

\noindent (b) $\delta_D(y) > R_2$ and $|x-y| \le R_2/L_1$: In this
case, $|X_{\tau_{B(y, R_2)}}-y| \,\ge \, R_2 \,\ge \,L_1 |x-y|$ and,
by the monotonicity of $G$ and Lemma \ref{G_G}, we get
\begin{eqnarray*}
\lefteqn{G_D(x,y) \,\ge\, G_{B(y, R_2)}(x,y)
\,
=\, G(x,y) - \E_x\left[     G( X_{\tau_{B(y, R_2)}}, y) \right]}\\
&\ge & G(x-y) - G(L_1 (x-y)) \,\ge\, c_1
\frac1{|x-y|^{d-\alpha}\ell(|x-y|^{-2})}, \qquad \forall |x-y| \le
\frac{R_2}{L_1}\, .
\end{eqnarray*}

\noindent (c) $\delta_D(y) > R_2$ and $|x-y| > R_2/L_1$: Choose a point $w \in
\partial B(y , R_2/L_1)$. Then from the argument in (b), we get
$$
G_D(w,y)  \ge G(w,y) - \E_w\left[     G( X_{\tau_{B(y, R_2)}}, y)
\right] \,\ge\, c_1
\frac1{(R_2/L_1)^{d-\alpha}\ell((R_2/L_1)^{-2})}.
$$
Since $D$ is bounded and $G_D(\,\cdot\,, y)$ is harmonic with respect to
$X$ in $B(x, R_2/(2L_1)) \cup B(w, R_2/(2L_1))$,  by Theorem \ref{HP2} we have
\begin{eqnarray*}
\lefteqn{G_D(x,y) \,\ge\, c_2 G_D(w,y)
\,\ge\, c_3 \frac1{(R_2/L_1)^{d-\alpha}\ell((R_2/L_1)^{-2})}}\\
&\ge & \frac{c_3}{ \ell((R_2/L_1)^{-2})} \left(\inf_{\frac{R_2}{L_1}
\le s \le M} \ell(s^{-2})\right)
\frac1{|x-y|^{d-\alpha}\ell(|x-y|^{-2})}, \qquad \forall |x-y| >
\frac{R_2}{L_1}\, .
\end{eqnarray*}
\qed

\begin{lemma}\label{G:g3}
For every $L >0$ and bounded open set $D$, there exists $C_{11}>0$ such
that for every $|x-y| \le L ( \delta_D(x) \wedge \delta_D(y)) $,
\begin{equation}\label{e:g3}
G_D(x,y) \ge C_{11} \frac1{|x-y|^{d-\alpha}\ell(|x-y|^{-2})}.
\end{equation}
\end{lemma}

\pf Without loss of generality, we assume $\delta_D(x) \le
\delta_D(y)$. Moreover, by Proposition \ref{GI} we can assume that
$L > 1/L_1$ and we only need to show \eqref{e:g3} for $\frac1{L_1}
\delta_D(x) \le |x-y| \le L \delta_D(x)$.

Choose a point $w \in \partial B(x , \delta_D(x)/L_1)$. Then by
Proposition \ref{GI}, we get
$$
G_D(x,w) \,\ge\, c_1
\frac1{(\delta_D(x)/L_1)^{d-\alpha}\ell((\delta_D(x)/L_1)^{-2})}.
$$
Since $|y-w| \le |x-y|+|x-w| \le (L+1) \delta_D(x)$ and
$G_D(x,\,\cdot\,)=G_D(\,\cdot\,, x)$ is harmonic with respect to $X$ in
$B(y, \delta_D(x)/L_1) \cup B(w, \delta_D(x)/L_1)$,
by Theorem \ref{HP2} we have
\begin{eqnarray}
\lefteqn{G_D(x,y) \ge c_2 G_D(x,w)
\ge c_3 \frac1{(\delta_D(x)/L_1)^{d-\alpha}\ell((\delta_D(x)/L_1)^{-2})}} \nonumber\\
& \ge & c_4 \left(  \frac{   \ell(|x-y|^{-2})} {\ell((\delta_D(x)/L_1
)^{-2})}  \right) \frac1{|x-y|^{d-\alpha}\ell(|x-y|^{-2})}.
\label{e:ggg3}
\end{eqnarray}
By the uniform convergence theorem (\cite[Theorem 1.2.1]{BGT}), we
can choose a small $r_1=r_1(\ell, L)>0$ such that
\begin{equation}\label{e:gg3}
\inf_{\lambda \in [1, 2L]} \frac{\ell((\lambda r)^{-2})}{
\ell(r^{-2})} \, \ge \, \frac12, \qquad \forall r \le r_1\, .
\end{equation}
If $\frac1{L_1} \delta_D(x) \le |x-y| \le L \delta_D(x)
\le r_1$, by \eqref{e:ggg3}--\eqref{e:gg3},
$$
G_D(x,y) \,\ge\, c_4 \left(  \inf_{r\le r_1 , \lambda \in [1, 2L]}
\frac{\ell((\lambda r)^{-2})}{ \ell(r^{-2})} \right)
\frac1{|x-y|^{d-\alpha} \ell(|x-y|^{-2})} \,\ge\, \frac{c_4}2
\frac1{|x-y|^{d-\alpha}\ell(|x-y|^{-2})} .
$$

On the other hand, if  $\frac1{L_1} \delta_D(x) \le |x-y| \le
L \delta_D(x)$ and $\delta_D(x) \ge r_1/L_1$, then
$|x-y| \ge\frac{r_1}{LL_1}.$ Thus from \eqref{e:ggg3}, we see that
 $$
G_D(x,y) \,\ge\, c_4 \left(  \inf_{r \in [\frac{r_1}{LL_1}, M]}
\ell(r^{-2}) \right)
 \left(  \inf_{r \in [\frac{r_1}{LL_1}, M]} \ell(r^{-2})^{-1} \right)
 \frac1{|x-y|^{d-\alpha}\ell(|x-y|^{-2})} $$
where $M=\mathrm{diam}(D)$. \qed

For the remainder of this section, we assume that $D$ is a bounded
$\kappa$-fat open set with characteristics $(R_1, \kappa)$.
Without loss of generality we may assume that $R_1\le 1/4$.
We recall that  for each $Q \in \partial D$ and $r \in (0, R_1)$,
$A_r(Q)$ denotes a point in $D \cap B(Q,r)$ satisfying $B(A_r(Q),\kappa r)
\subset D \cap B(Q,r)$.
Recall also that $G_D(\cdot, z)$ is
regular harmonic in $D\setminus \overline{B(y,\varepsilon)}$ for
every $\varepsilon >0$ and vanishes outside $D$. From Theorem
\ref{BHP}, we get the following boundary Harnack principle for the
Green function of
$X$ in $D$.

\begin{thm}\label{BHP2}
There exists a constant $C_{12}>1$ such that for any $Q \in \partial D$,
$r \in (0,R_1]$ and $z,w \in D \setminus B(Q,2r)$, we have
$$
C_{12}^{-1}\, \frac{G_D(z, A_r(Q))}{G_D(w, A_r(Q))} \le
\frac{G_D(z,x)}{G_D(w, x)}
\le C_{12}\, \frac{G_D(z,A_r(Q))} {G_D(w, A_r(Q))} ,
\quad x\in D\cap B\left(Q,\frac{r}{2}\right)\, .
$$
\end{thm}

Using the uniform convergence theorem (\cite[Theorem 1.2.1]{BGT}),
we further choose $R_3\le R_1$  such that
\begin{equation}\label{e:rsmall}
\frac12 \, \le \, \min_{ \frac16  \le \lambda
\le 2 \kappa^{-1}} \frac{ \ell((\lambda r)^{-2})}{
\ell(r^{-2})}\,\le\, \max_{ \frac16  \le \lambda \le 2 \kappa^{-1}}
\frac{ \ell((\lambda r)^{-2})}{ \ell(r^{-2})} \, \le\, 2, \quad
\mbox{ if } r\le R_3\, .
\end{equation}

Let us fix $z_0 \in D$ with $\kappa R_3  < \delta_D(z_0) < R_3$
and let $\eps_1:=  \kappa R_3/24$. For $x,y \in D$, we define
$r(x,y): = \delta_D(x) \vee \delta_D(y)\vee |x-y|$ and
\begin{equation} \label{d:gz1}
\BB(x,y):=
\begin{cases}
\left\{ A \in D:\, \delta_D(A) > \frac{\kappa}{2}r(x,y), \,
|x-A|\vee
|y-A| < 5 r(x,y)  \right\}& \text{ if } r(x,y) <\eps_1 \\
\{z_0 \}& \text{ if } r(x,y) \ge \eps_1 .
\end{cases}
\end{equation}
Note that for every $(x,y) \in D \times D$ with $r(x,y) <\eps_1$
\begin{equation}\label{e:ar}
\frac16 \delta_D(A) \,\le\,\delta_D(x) \vee
\delta_D (y) \vee |x-y|\, \le\, 2 \kappa^{-1}  \delta_D(A), \qquad A
\in \BB(x,y).
\end{equation}
So by \eqref{e:rsmall}, if $r(x,y) <\eps_1 $
\begin{equation}\label{e:areq}
\frac12 \, \le \, \frac{\ell((\delta_D(A))^{-2})}
{\ell((r(x,y))^{-2})} \,\le\, 2, \qquad A \in \BB(x,y).
\end{equation}
Let
\begin{equation}\label{e:C_1}
C_{13}:=C_{10} 2^{d-\alpha} \delta_D(z_0)^{-d+\alpha}
\left({\sup_{\delta_D(z_0)/2\le r\le M}\ell(r^{-2})^{-1}}\right)\, .
\end{equation}
It follows from Proposition \ref{GI} that $G_D(\cdot, z_0)$ is
bounded above by $C_{13}$ on $ D \setminus B(z_0, \delta_D(z_0)/2)$.
Now we define
\begin{equation}\label{d:gz0} g(x ):=  G_D(x, z_0) \wedge C_{13}.
\end{equation}
Note that if $\delta_D(z) \le 6 \eps_1$, then $|z-z_0| \ge
\delta_D(z_0) - 6 \eps_1 \ge \delta_D(z_0) /2$ since
$6\eps_1<\delta_D(z_0)/4$, and therefore $g(z )= G_D(z, z_0)$.

The two lemmas below follow immediately from Theorem \ref{HP2}.
\medskip

\begin{lemma}\label{G:g1}
There exists $C_{14}=C_{14}(\kappa)>1$ such that for all $x, y \in D$,
$C_{14}^{-1} g(A_1) \le g(A_2)\le C_{14} g(A_1)$ for all $A_1, A_2 \in \BB(x,y)$.
\end{lemma}

\begin{lemma}\label{G:g2}
There exists $C_{15}>1$ such that for every $x \in \{ y \in D:\,
\delta_D(y) \ge \kappa^3\eps_1 /64\} $, $C_{15}^{-1}  \le g(x) \le C_{15}$.
\end{lemma}

\medskip
We are now ready to prove Theorem \ref{t:Gest}.
\medskip

\noindent {\bf Proof of  Theorem \ref{t:Gest}}.
Since the proof is an adaptation of the proofs of
\cite[Proposition 6]{B3} and \cite[Theorem 2.4]{H}, we only give the
proof when $\delta_D(x) \le \delta_D(y)\le \frac{\kappa}{4} |x-y|$.

In this case, we have $r(x,y)=|x-y|$. Let $r:= \frac12 (|x-y| \wedge
\eps_1)$. Choose $Q_x, Q_y \in \partial D$ with $|Q_x-x|
=\delta_D(x)$ and $|Q_y-y| =\delta_D(y)$. Pick points $x_1=A_{\kappa
r/2}(Q_x)$ and $y_1=A_{\kappa r/2}(Q_y)$ so that $x, x_1 \in B(Q_x,
\kappa r/2)$ and $y, y_1 \in B(Q_y, \kappa r/2)$. Then one can
easily check that $|z_0-Q_x| \ge r$ and $|y-Q_x| \ge r$. So by
Theorem \ref{BHP2}, we have
$$
c_1^{-1}\,  \frac{G_D(x_1,y)}{g(x_1)}\,\le\, \frac{G_D(x,y)}{g(x)}
\,\le\, c_1 \frac{G_D(x_1,y)}{g(x_1)}
$$
for some $c_1>1$. On the other hand, since  $|z_0-Q_y| \ge r$ and
$|x_1-Q_y| \ge r$, applying Theorem \ref{BHP2} again,
$$
c_1^{-1}\,\frac{G_D(x_1,y_1)}{g(y_1)} \,\le\,
\frac{G_D(x_1,y)}{g(y)} \,\le\, c_1 \frac{G_D(x_1,y_1)}{g(y_1)}.
$$
Putting the two  inequalities above together we get
$$
c_1^{-2}\,\frac{G_D(x_1,y_1)}{g(x_1)g(y_1)}  \,\le\,
\frac{G_D(x,y)}{g(x)g(y)} \,\le\, c_1^2\frac{G_D(x_1,y_1)}
{g(x_1)g(y_1)}.
$$
Moreover, $\frac13|x-y| < |x_1-y_1| < 2 |x-y| $ and $|x_1 -y_1| \le
\frac{64}{\kappa^3\eps_1} (\delta_D(x_1) \wedge \delta_D(y_1))$.
Thus by Lemma \ref{G:g3}, we have
\begin{equation}\label{e:GE1}
\frac{2^{-d+\alpha}c_2^{-1} c_1^{-2}}{g(x_1)g(y_1)}
\frac1{|x-y|^{d-\alpha}\ell(|x_1-y_1|^{-2})}    \,\le\,
\frac{G_D(x,y)}{g(x)g(y)} \,\le\, \frac{3^{d-\alpha}c_2
c_1^2}{g(x_1)g(y_1)}\frac1{|x-y|^{d-\alpha}\ell(|x_1-y_1|^{-2})}
\end{equation}
for some $c_2>1$.

If $r=\eps_1/2$, then $r(x,y)=|x-y| \ge \eps_1$. Thus
$g(A)=g(z_0)=C_{13}$ and $\delta_D(x_1)\wedge \delta_D(y_1) \ge
\kappa r/2 = \kappa \eps_1/4$. So by Lemma \ref{G:g2},
\begin{equation}\label{e:GE2}
C_{13}^{-2} c_3^{-2} \le\frac{g(A)^2}
{g(x_1)g(y_1)} \le C_{13}^2 c_3^2
\end{equation}
for some $c_3>1$.

If $r<\eps_1/2$, then $r(x,y)=|x-y| < \eps_1$ and $r=\frac12
r(x,y)$. Hence $\delta_D(x_1), \delta_D(y_1) \ge \kappa r/2 = \kappa
r(x,y) /4$. Moreover, $|x_1-A|, |y_1-A| \ge 6 r(x,y)$. So by
applying Theorem \ref{T:Har} to $g$ twice,
\begin{equation}\label{e:GE3}
c^{-1}_4 \,\le\,\frac{g(A)}{g(x_1)} \,\le\,c_4 \quad \mbox{and}
\quad c^{-1}_4 \,\le\,\frac{g(A)}{g(y_1)} \,\le\,c_4
\end{equation}
for some constant $c_4=c_4(D)>0$. Combining
\eqref{e:GE1}-\eqref{e:GE3}, we get
$$
c_5^{-1}\,\frac{g(x) g(y)}{g(A)^2|x-y|^{d-\alpha}\ell(|x_1-y_1|^{-2})}\le G_D(x,y) \le c_5\,\frac{g(x)
g(y)}{g(A)^2|x-y|^{d-\alpha}\ell(|x_1-y_1|^{-2})}, \quad  A \in \BB(x,y)\, .
$$
If  $\frac13|x-y| < |x_1-y_1| < 2 |x-y| \le R_3$, by \eqref{e:rsmall},
$\frac12  \le  \frac{\ell(|x-y|^{-2})} {\ell(|x_1-y_1|^{-2})}\le 2 $.
On the other hand, if $\frac13|x-y| < |x_1-y_1| < 2 |x-y| $ and $ 2
|x-y| > R_3$, we have that $ \frac{R_3}2 \le |x-y|\le M $ and $\frac{R_3}{6}\le |x_1-y_1| \le M$.
Thus from the local boundedness  of $\ell$ on $(0, \infty)$, we see that
 $$
c_6^{-1} \, \le\, \frac{\ell(|x-y|^{-2})}
{\ell(|x_1-y_2|^{-2})}\,\le\,   c_6
$$
for some $c_6>1$. \qed

\section{Explicit Green function estimates on bounded $C^{1,1}$-open sets}\label{s:4}

In this section we refine the estimates from Theorem \ref{t:Gest} in
the case of  bounded $C^{1,1}$ open sets.

Recall that $X=(X_t:\, t\ge 0)$ is the $d$-dimensional subordinate
Brownian motion defined by $X_t=W_{S_t}$ where $W=(W^1,\dots, W^d)$
is a $d$-dimensional Brownian motion and $S=(S_t:\, t\ge 0)$ an
independent subordinator with the Laplace exponent $\phi$
which is a complete Bernstein function satisfying assumption {\bf (H)}
 for $\alpha\in (0, 2\wedge d)$ and the additional
assumption  \eqref{e:ass4trans}
when $d\le 2$.
Let $Z=(Z_t:\, t\ge 0)$ be
the one-dimensional subordinate Brownian motion defined as $Z_t:=W^d_{S_t}$.
Recall that the potential measure of the ladder height process of
$Z$ is denoted by $V$ and its density by $v$.
We also use $V$ to denote the renewal function of the ladder height process of $Z$.
In Corollary \ref{c:vsm} we have established that both $V$ and
$v$ are $C^{\infty}$ function.
Recall the notation $\bR^d_+:=\{
x=(x_1, \dots, x_{d-1}, x_d):=(\tilde{x}, x_d)  \in \bR^d: x_d > 0
\}$ for the half-space. The next result, which follows from
\cite[Theorem 2]{Sil}, is very important in this paper.

Define $w(x):=V((x_d)^+)$.

\begin{thm}\label{t:Sil}
The function $w$ is harmonic in $\R^d_+$ with respect to $X$
and, for any $r>0$, regular harmonic in
$\R^{d-1}\times (0, r)$ (in $(0,r)$ when $d=1$)
  with respect to $X$.
\end{thm}
\pf
Since $Z_t:=W^d_{S_t}$ has a transition density, it
satisfies the condition ACC in  \cite{Sil}, namely the resolvent
kernels are absolutely continuous. The assumption in  \cite{Sil}
that $0$ is regular for $(0,\infty)$ is also satisfied since $Z$ is
symmetric and has infinite L\'evy measure.  Indeed, if $0$ were
irregular for $(0,\infty)$, it would be, by symmetry, irregular
for $(-\infty,0)$ as well. But then $Z$ would be a compound
Poisson process which contradicts the fact that it has infinite L\'evy measure.
Further, again by symmetry of $Z$, the notions
of coharmonic and harmonic functions coincide.
Let $Z^{(0,\infty)}$ (respectively $X^{\R^d_+}$) denote the
process $Z$ killed upon exiting $(0,\infty)$ (respectively $X$
killed upon exiting $\R^d_+$).
By \cite[Theorem 2]{Sil}, the renewal function $V$ of the ladder
height process of $Z$ is invariant for $Z^{(0, \infty)}$.
Thus $w$ is invariant for $X^{\bR^d_+}$.
Being invariant for $X^{\R^d_+}$, $w$ is also harmonic for $X^{\R^d_+}$, and consequently,
harmonic in $\R^d_+$ with respect to $X$.
We show now that $w$ is regular harmonic for $X$ in $\R^{d-1}\times (0, r)$ for any $r>0$.
First note that since $V$ is continuous at zero and $V(0)=0$, it follows that
\begin{equation}\label{e:w-to-0}
\lim_{x_d\to 0}w(x)=\lim_{x_d\to 0}w(\tilde{x},x_d)=\lim_{x_d\to 0} V(x_d)=0\, .
\end{equation}
Thus,  by harmonicity of $w$ and \eqref{e:w-to-0}
$$
w(x)=w(\tilde{x}, x_d)=\lim_{\eps \to 0} \E_{x}\left[w\left(
X_{\tau_{ \R^{d-1}\times (\varepsilon,r)}}\right)\right]=
\E_{x}\left[w\left(X_{\tau_{ \R^{d-1}\times (0,r)}}\right)\right]\,
, \quad x_d>0\, .
$$
\qed

\begin{prop}
\label{c:cforI}
For all positive constants
$r_0$ and $L$, we have
$$
\sup_{x \in \R^d:\, 0<x_d <L} \int_{B(x, r_0)^c \cap \bR^d_+}
w(y) j(|x-y|)\, dy < \infty\, .
$$
\end{prop}
\pf
By Theorem \ref{t:Sil} and \eqref{newls},
for every $x \in \R^d_+$,
\begin{eqnarray}
w(x)
&\ge& \E_x\left[w(X_{\tau_{ B(x, r_0/2)\cap \bR^d_+}}):
X_{\tau_{ B(x, r_0/2)\cap \bR^d_+}} \in  {B(x, r_0)}^c \cap
\bR^d_+  \right]\nn\\
&=& \int_{{B(x, r_0)}^c \cap \bR^d_+}\int_{B(x, r_0/2)\cap \bR^d_+}
G_{B(x, r_0/2)\cap \bR^d_+} (x,z) j(|z-y|)w(y)\, dz\, dy\, .
\label{e:fdfsaff}
\end{eqnarray}
Since $|z-y| \le |x-z|+|x-y| \le r_0 +|x-y| \le 2|x-y|$ for $(z,y)
\in B(x, r_0/2) \times B(x, r_0)^c$, using \eqref{H:1} and
\eqref{H:2}, we have $j(|z-y|) \ge c_1 j(|x-y|)$.
Thus, combining this with \eqref{e:fdfsaff}, we obtain that
$$
 \sup_{x \in \R^d:\, 0<x_d <L}
\int_{B(x, r_0)^c \cap \bR^d_+} w(y) j(|x-y|) dy  \le
c_1^{-1}\sup_{\wt x=0,  0<x_d <L}  \frac{w(x)}{\E_{x}[\tau_{B(x, r_0/2)\cap \bR^d_+}]}.
$$
We claim that the supremum on the right-hand side above is finite.
Clearly, if $L>x_d \ge r_0/(64)$ and $\wt x=0$,
$$
\frac{w(x)}{\E_{x}[\tau_{B(x, r_0/2)\cap \bR^d_+}]} \le
\frac{V(L)}{\E_{0}[\tau_{B(0, r_0/(64))}]}.
$$

Suppose $x_d  < r_0/(64)$  and $\wt x=0$.
Let $U:=B((\wt 0, 16r_0), r_0)$. By the L\'evy system, we have
$$
\P_x\left(X_{\tau_{B(0, r_0/2)\cap \bR^d_+}} \in U\right) = \int_U
\int_{B(0, r_0/2)\cap \bR^d_+} G_{B(0, r_0/2)\cap \bR^d_+}(x,z)
 j(|z-y|)dz dy \le c_2 \E_{x} [\tau_{B(0, r_0/2)\cap \bR^d_+}].
$$
Thus, by the above and the boundary Harnack principle (Theorem \ref{BHP}),
$$
\frac{w(x)}{\E_{x} [\tau_{B(x, r_0/2)\cap \bR^d_+}]}\le
 c_2
\frac{w(x)}{\P_x(X_{\tau_{B(0, r_0/2)\cap \bR^d_+}} \in U)}
\le
 c_3
\frac{w(x_1)}{\P_{x_1}(X_{\tau_{B(0, r_0/2)\cap \bR^d_+} }\in U)}
\le c_4V(r_0/(16))
$$
where $x_1=(\wt 0, r_0/(16))$. We have thus proved the claim.
\qed

We now define the
operator
($\AA$, $\mathfrak{D}(\AA)$)
by the following formula:
\begin{eqnarray}
 \AA  f(x)&:=& \mathrm{P.V.} \int_{\R^d}
\left(f(y)-f(x)\right)j(|y-x|)\, dy:=\lim_{\eps \downarrow 0}
\int_{\{y\in \bR^d: |x-y| > \eps\}}
\left(f(y)-f(x)\right)j(|y-x|)\, dy\, \nn\\
\mathfrak{D}(\AA)&:=&\left\{f:\R^d\to \R: \lim_{\eps \downarrow 0}
\int_{\{y\in \bR^d: |x-y| > \eps\}}
\left(f(y)-f(x)\right)j(|y-x|)\, dy \text{ exists and is finite } \right\}.
\label{generator}
\end{eqnarray}
 It is well known that
$C^2_0\subset \mathfrak{D}(\AA)$, where  $C^2_0$ is the collection of $C^2$ functions in
$\RR^d$ vanishing at infinity, and
that by the rotational symmetry of $X$, $\AA$ restricted to $C^2_0$
coincides with the infinitesimal generator of the process $X$ (e.g.
\cite[Theorem 31.5]{Sa}).

\begin{thm}\label{c:Aw=0}
$\AA w(x)$ is well defined and $\AA w(x)=0$ for all $x \in \bR^d_+$.
\end{thm}
\pf
We first note that it follows from
Proposition \ref{c:cforI} and the fact that $j$ is a L\'evy density
that for any $L>0$
\begin{eqnarray}
\lefteqn{\sup_{x \in \R^d:\ 0<x_d <L}
\left|\int_
{\{y  \in \R^d:|y-x|>1\}}(w(y)-w(x)) j(|y-x|)dy \right|\nn}\\
&\le&
\sup_{x \in \R^d:\ 0<x_d <L}  \int_{\{y  \in \R^d:|y-x|>1\}}w(y)j(|y-x|)\, dy
+V(L)\int_{\{y  \in \R^d:|y|>1\}}j(|y|)dy
<\infty\, . \label{e:dedwer}
\end{eqnarray}
Hence, for every $\varepsilon \in (0,1/2)$
$$
\AA_\eps w(x):=\int_{\{y  \in \R^d:|y-x|>\varepsilon\}}(w(y)-w(x))j(|y-x|)dy
$$
is well defined.
Note that since $w(x)=V((x_d)^+)$ and $V$ is smooth in $(0,\infty)$
by Corollary \ref{c:vsm}, it holds that $w$ is smooth in $\R^d_+$.  Hence,
$$
\AA_\eps w(x)=\int_{\{y  \in \R^d:|y-x|>\varepsilon\}}\left(w(y)-w(x)-{\bf 1}_
{\{|y-x|<1\}}
(y-x)
\cdot\nabla w(x)\right)j(|y-x|) dy.
$$
Moreover, by the smoothness of $w$,
$$
x \mapsto \int_
{\{y  \in \R^d:|y-x| \le \varepsilon\}}
\left(w(y)-w(x)-(y-x)
\cdot\nabla w(x)\right)j(|y-x|)\, dy
$$
converges to 0 locally uniformly in $\R^d_+$ as $\eps \to 0$ .
Combining this with \eqref{e:dedwer}, we see that $\AA w$ is
well defined in $\R^d_+$ and
$\AA_\eps w(x)$
converges to
$$\AA w(x)=\int_{\RR^d}\left(w(y)-w(x)-{\bf 1}_
{\{|y-x|<1\}}
(y-x)
\cdot\nabla w(x)\right)j(|y-x|) dy$$ locally uniformly
in ${\R^d_+}$ as $\varepsilon\to 0$.

Moreover, for every $x \in \RR^d_+$, $z \in B(x, (\eps
\wedge x_d)/2)$, and $y\in B(z,\varepsilon)^c$ it holds that
$\frac12 |y-z| \le |y-x| \le \frac32 |y-z|.$
So,  using \eqref{H:1},
\begin{align*}
\lefteqn{{\bf 1}_
{\{|y-z|>\varepsilon\}}
\left|(w(y)-w(z)-{\bf 1}_
{\{|y-z|<1\}}
(y-z) \cdot\nabla w(z))\right|j(|y-z|) }\\
\le& c \left(\sup_{\eps/2 < s< x_d+2} V''(s)\right) |y-x|^2 {\bf 1}_
{\{\eps/2<|y-x|<2\}}\, j(|y-x|/2)\\
&
\quad + (w(y) +V(x_d+1)) {\bf 1}_
{\{|y-x|>1/2\}}
\, j(|y-x|/2).
\end{align*}
It follows from Proposition \ref{c:cforI} and the fact that
$j$ is a L\'evy density, by using the dominated convergence theorem,
that $x \to \AA_\eps w(x)$ is continuous for each $\eps$. Therefore,
by this and the local uniform convergence of $\AA_\eps w$, the
function $\AA w(x)$ is continuous in ${\R^d_+}$.

Suppose that $U_1$ and $U_2$ are relatively compact open subsets of $\RR^d_+$ such that
$\overline{U_1} \subset U_2
\subset \overline{U_2} \subset \RR^d_+$. Let $r_0:=\mathrm{dist}(U_1, U_2^c)>0$.
Then, by Proposition \ref{c:cforI}
 \begin{eqnarray}
\int_{U_1} \int_{U_2^c} w(y) j(|x-y|)dydx
&\le& |U_1| \sup_{x \in U_1} \int_{U_2^c} w(y) j(|x-y|)dy \nn\\
&\le& |U_1| \sup_{x \in U_1} \int_{B(x, r_0)^c} w(y) j(|x-y|)dy <\infty\, . \label{e:(2.4)}
\end{eqnarray}
By harmonicity of $w$, clearly $w( X_{\tau_{U_1}}) \in L^1(\P_x)$ and
$$\sup_{x \in U_1} \E_x\left[{\bf 1}_{U_2^c} (X_{\tau_{U_1}}) w( X_{\tau_{U_1}})\right]
\le \sup_{x \in U_1} \E_x\left[ w( X_{\tau_{U_1}})\right] =\sup_{x \in U_1} w(x) <\infty.
$$
The last two displays show that the conditions \cite[(2.4), (2.6)]{C} are true.
Thus, by \cite[Lemma 2.3, Theorem 2.11(ii)]{C}, we have that for any $f\in C^2_c(\R^d_+)$,
\bee\label{e:C2.11}
0=\int_{\R^d}\int_{\R^d}(w(y)-w(x))(f(y)-f(x))j(|y-x|)\, dx\, dy.
\eee
For $f\in C^2_c(\R^d_+)$ with supp$(f) \subset \overline{U_1} \subset U_2
\subset \overline{U_2} \subset \RR^d_+$,
\begin{align*}
&\int_{\R^d}\int_{\R^d}|w(y)-w(x)||f(y)-f(x)|j(|y-x|) dxdy\\
=&\int_{U_2}\int_{U_2}|w(y)-w(x)||f(y)-f(x)|j(|y-x|) dxdy+2\int_{U_1}\int_{U_2^c}|w(y)-w(x)||f(x)|j(|y-x|) dxdy\\
\le&c_1\int_{U_2 \times U_2}
|y-x|^2 j(|y-x|) dxdy+2\|f\|_\infty|U_1|  \left(\sup_{x \in U_1} w(x)\right) \int_{U_2^c}j(|y-x|) dy\\
&+2\|f\|_\infty\int_{U_1} \int_{U_2^c} w(y) j(|x-y|)dydx
\end{align*}
is finite by \eqref{e:(2.4)} and the fact that $j$ is a L\'evy density.
Thus by \eqref{e:C2.11}, Fubini's theorem and the dominated convergence theorem,
we have
for any $f\in C^2_c(\R^d_+)$,
\begin{eqnarray*}
&&0=\lim_{\varepsilon\downarrow0}\int_{\{(x, y)\in {\R^d}\times {\R^d},\
|y-x|>\varepsilon\}}(w(y)-w(x))(f(y)-f(x))j(|y-x|)\, dx\, dy\\
&&=
-2\lim_{\varepsilon\downarrow0}\int_{\R^d_+} f(x)
\left(\int_{\{y \in \RR^d:|y-x|>\varepsilon\}}(w(y)-w(x))j(|y-x|) dy \right) dx\,=\,
-2\,\int_{\R^d_+} f(x)\AA w(x)\, dx,
\end{eqnarray*}
where we have used the fact $\AA_\eps w \to \AA w$ converges uniformly on the support of $f$.
Hence, by the continuity of $\AA w$, we have $\AA w(x)=0$ in ${\R^d_+}$.
\qed

For $x\in \bR^d$, let $\delta_{\partial D}(x)$ denote the Euclidean
distance between $x$ and $\partial D$.
It is well known that any
$C^{1, 1}$ open set $D$ with characteristics $(R, \Lambda)$ satisfies both the {\it uniform interior
ball condition} and the {\it uniform exterior ball condition} with
the radius $r_1$: there exists $r_1< R$ such that for every $x\in D$
with $\delta_{D}(x)< r_1$ and $y\in \bR^d \setminus \overline D$
with $\delta_{D}(y)<r_1$, there are $z_x, z_y\in \partial D$ so that
$|x-z_x|=\delta_{\partial D}(x)$, $|y-z_y|=\delta_{D}(y)$ and that
$B(x_0, r_1)\subset D$ and $B(y_0, r_1)\subset \bR^d \setminus
\overline D$ for $x_0=z_x+r_1(x-z_x)/|x-z_x|$ and
$y_0=z_y+r_1(y-z_y)/|y-z_y|$.

In the remainder of this section, we assume $D$ is a $C^{1,1}$ open
set with characteristics $(R, \Lambda)$
and $D$ satisfies the  uniform
interior ball condition and the uniform exterior ball condition with
the radius $R$ (by choosing $R$ smaller if necessary).

\begin{lemma}\label{L:Main}
Fix $Q \in \partial D$ and let
$$
h(y):=V\left(\delta_D (y)\right){\bf 1}_{D\cap B(Q, R)}(y).
$$
There exist
$C_{16}=C_{16}(\alpha, \Lambda, R, \ell)>0$ and $R_4 \le R/4$ independent of
the point $Q \in \partial D$ such that $\AA h$ is well defined in $D\cap B(Q, R_4)$ and
\bee\label{e:h3}
|\AA h(x)|\le C_{16} \quad \text{ for all } x \in D\cap B(Q, R_4)\, .
\eee
\end{lemma}

\pf
We first note that when $d=1$, the lemma follows from Proposition
\ref{c:cforI} and Theorem \ref{c:Aw=0}. In fact, suppose that  $d=1$ and
$x \in D\cap B(Q, R/2)$.
Without loss of generality we may assume that $Q$ is the origin and
$D\cap B(Q, R)=(0, R)$ (due to uniform exterior ball condition).
Since $h(y)=w(y)$ for $y \in D\cap B(Q, R)=(0, R)$, we have

\begin{eqnarray*}
 \AA (h-w)(x)
&=& \lim_{\eps \downarrow 0}
\int_{\{y\in \bR^1: |x-y| > \eps\}}
(h-w) (y) j(|y-x|)\, dy\\
&=& \lim_{\eps \downarrow 0}
\int_{\{y \in (0, R)^c: |x-y| > \eps \}}
(h-w) (y)j(|y-x|)\, dy\\
&=& -\lim_{\eps \downarrow 0}
\int_{\{y \in (0, R)^c: |x-y| > \eps \}}
w (y) j(|y-x|)\, dy\\
&=& -\lim_{\eps \downarrow 0}
\int_{\{y \ge R: |x-y| > \eps \}}
w (y) j(|y-x|)\, dy\\
\end{eqnarray*}
Since $0<x<R/2$ and $y \ge R$, we have $|x-y| >R/2$, and thus
  \begin{eqnarray*}
|\AA (h-w)(x)|   \le \int_{\{y \ge R,  |x-y| > R/2 \}}
w (y) j(|y-x|)\, dy.
\end{eqnarray*}
Therefore, by using Theorem \ref{c:Aw=0} (which gives $\AA h(x)= \AA (h-w)(x)$),
Proposition \ref{c:cforI} and the above display, we conclude that
$$|\AA h(x)|= |\AA (h-w)(x)| \le
\sup_{0<z<R/2}
\int_
{\{y  \in \R^1_+:|y-z|>R/2\}} w(y) j(|y-z|)dy  < \infty.
$$

Throughout the remainder of the proof,
$d \ge 2$
and $q$ is a fixed positive constant such that
$$
0<q<\frac{\alpha \wedge(2-\alpha)}{20}
\quad\text{and} \quad \alpha+2q-1\neq 0.$$

Since $\ell$ is slowly varying at $\infty$, by Potter's Theorem
(\cite[Theorem 1.5.6 (1)]{BGT}), we can find a small $R_4<1 \wedge
(R/4)$ such that
 for every $r \le 2 R_4^2$
\begin{eqnarray}
\frac{\ell( r^{-2})}{(\ell( (2 R^{-1})^{-2} r^{-4}))^{1/2}}&\le&
2\,\frac{\ell( (2 R_4^2)^{-2})}{(\ell( (2 R^{-1})^{-2} (2 R_4^2)^{-4}))^{1/2}}
(2 R_4^2)^{1/2} r^{-1/2}   \le c_1  r^{-1/2},\label{e:upas1}\\
\ell( r^{-1}) &\le& 2\,\ell( (2 R_4^2)^{-1}) (2 R_4^2)^{q} r^{-q}
\le c_1  r^{-q},\label{e:upas2}\\
\ell( r^{-2})^{-1/2} &\le& 2\,  \ell( (2 R_4^2)^{-2})^{1/2} (2
R_4^2)^{q}r^{-q}\le c_1  r^{-q}.\label{e:upas3}
\end{eqnarray}

In the remainder of this proof, we fix $x \in D\cap B(Q, R_4)$ and
$x_0\in\partial D$ satisfying $\delta_D(x)=|x-x_0|$.  We also fix
the $C^{1, 1}$ function $\psi$ and the coordinate system $CS=CS_{x_0}$
in the definition of $C^{1, 1}$ open set
so that $x=( 0, x_d)$ with $0<x_d <R_4$ and
$B(x_0, R)\cap D=\{ y=(\wt y, \, y_d) \in B(0, R) \mbox{ in } CS :
y_d > \psi (\wt y) \}.$
Let
$$
\psi_1(\wt y):=R -\sqrt{ R^2-|\wt y|^2}
\quad \text{and}\quad \psi_2(\wt y):=-R +\sqrt{ R^2-|\wt
y|^2}.
$$
Due to the  uniform interior ball condition and the uniform
exterior ball condition with the radius $R$, we have
 \bee\label{e:phi012}
\psi_2(\wt y)  \le \psi (\wt y) \le \psi_1(\wt y) \quad \text{for
every } y \in D\cap B(x, R_4).
 \eee
Define
$H^+:=
\left\{y=(\wt y, \, y_d) \in CS:y_d>0 \right\}$ and
 let
$$
A:=\{y=(\widetilde{y},y_d) \in (D \cup H^+)\cap B(x, R_4):
\psi_2(\wt y)  \le y_d  \le \psi_1(\wt y)\},
$$
$$
E:=\{y=(\widetilde{y},y_d) \in
  B(x, R_4):  y_d  > \psi_1(\wt
y)\}.
$$
Note that, since $|y-Q| \le |y-x| +|x-Q| \le R/2$ for $y
\in B(x, R_4)$, we have
 \bee\label{e:BE} B(x,  R_4) \cap D \subset
B(Q,   R/2)\cap D \, .
 \eee

Let
$$
h_{x}(y):=V\left(\delta_{_{H^+}}(y)\right).
$$
Note that $h_{x}(x)=h(x)$. Moreover, since
$\delta_{_{H^+}}(y)=(y_d)^+$ in $CS$, by Theorem
\ref{c:Aw=0}
it follows that $\AA h_x$ is well defined in $H^+$ and
\bee\label{e:hz}
\AA  h_{x}(y)=0, \quad  \forall y\in H^+.
\eee
We show now that $\AA (h-h_x)(x)$ is well defined.
For each $\varepsilon >0$ we have that
\begin{align*}
&\bigg|\int_{\{y \in D \cup H^+: \, |y-x|>\varepsilon\}}
{(h(y)-h_{x}(y))}j(|y-x|)\ dy\bigg|
\nn\\
\leq & \int_{B(x,R_4)^c}
(h(y)+h_{x}(y))
j(|y-x|)dy +\int_{A} (h(y)+h_{x}(y))j(|y-x|)\ dy +\int_{E}
{|h(y)-h_{x}(y)|}j(|y-x|) dy\\
=:&I_1+I_2+I_3.
\end{align*}

We claim that
 \bee\label{e:I1234}
I_1+I_2+I_3 \leq C_{16}
 \eee
for some constant $C_{16}=C_{16}(\alpha,\Lambda, R, \ell)$.
This shows in particular that the limit
$$
\lim_{\varepsilon\downarrow 0}\int_{\{y \in D \cup H^+:
|y-x|>\varepsilon\}}{(h(y)-h_{x}(y))}j(|y-x|)\, dy
$$
exists and
hence $\AA(h-h_x)(x)$ is well defined, and
$|\AA(h-h_x)(x)|\le C_{16}$. By linearity
and \eqref{e:hz},
we get that $\AA h(x)$ is well defined and $|\AA h(x)|\le C_{16}$.
Therefore, it remains to prove \eqref{e:I1234}.

By the fact that $h(y)=0$ for $y \in B(Q, R)^c$,
\begin{align}
I_1
\le &\sup_{z \in \R^d:\ 0<z_d <R} \int_{B(z,  R_4)^c \cap H^+}    V(y_d) j(|z-y|)
dy+c_3\int_{B(0,  R_4)^c}j (|y|)dy
=:K_1+K_2.
\nonumber
\end{align}
$K_2$ is clearly finite since $J$ is the L\'evy density
of $X$ and $K_1$ is  finite by
Proposition \ref{c:cforI}.

For $y \in A$, since $V$ is increasing and
$(R - \sqrt{R^2-|\wt y|^2}) \le  R^{-1}|\wt y|^2$,
we see that
\begin{align}
 h_{x}(y)+h(y) \le 2V(\psi_1 ( \wt y) -\psi_2 ( \wt y)) \le2 V(2 R^{-1}|\wt y|^2)
 \label{e:dfe2}.
\end{align}
Using \eqref{e:dfe2},
\eqref{e:behofV} and Theorem \ref{t:Jorigin}, we have
\begin{align}
I_2\leq&\int_{ 0}^{ R_4} \int_{|\widetilde{y}|=r}{\bf 1}_A(y)
(h_{x}(y)+h(y))j(
(r^2+ |y_d- x_d|^2)^{1/2}
) \ m_{d-1}(d y)dr\nonumber\\
\le &
2\int_{ 0}^{  R_4}\int_{| \widetilde{y} |=r}{\bf 1}_A(y)V(2
R^{-1} r^2) j(r)\ m_{d-1}(d y)dr \nonumber\\
\le&
c_4\int_0^{ R_4}r^{-d}   \frac{\ell( r^{-2})}
{\ell( (2^{-2}R^{2}
r^{-4})^{1/2}}  m_{d-1}(\{y\in A:
|\widetilde{y}|=r\})dr \label{e:fgy6} \end{align}
where $m_{d-1}$ is the surface measure, that is, the $(d-1)$-dimensional Lebesgue measure.
Furthermore, since
$|\psi_2(\widetilde{y})-\psi_{1}(\widetilde{y})| \le 2
R^{-1}|\widetilde{y}|^{2}=2 R^{-1} r^2$ on $|\widetilde{y}|=r$,
we have for $r\le R_4$,
$m_{d-1}(\{ y:\, | \widetilde{y} |=r,\psi_2(\widetilde{y})<y_d<
\psi_{1}(\widetilde{y})\})\le
c_5r^d$
for some constant $c_5$.
Using the above inequality  and \eqref{e:upas1}, from \eqref{e:fgy6}
we get
\begin{align}
I_2\leq  c_6\int_0^{ R_4}\frac{\ell( r^{-2})}{(\ell( (2
R^{-1})^{-2} r^{-4}))^{1/2}}dr  \le c_7\int_0^{ R_4}r^{-1/2}dr
<\infty.\nonumber
\end{align}

For $I_3$, we consider two cases separately:
If
$0 <y_d=\delta_{_{H^+}}({y}) \le  \delta_D({y})$, since $v$ is
decreasing,
\begin{align}
h(y)-h_{x}(y) \le V(y_d+R^{-1}|\wt y|^2) -V(y_d) =
\int_{y_d}^{y_d+R^{-1}|\wt y|^2} v(z)dz \le R^{-1}|\wt y|^2
v(y_d).  \label{e:KG}
\end{align}
If $y_d=\delta_{_{H^+}}({y}) >  \delta_D({y})$ and $y \in E$, using the fact that
$\delta_D({y})$ is greater than or equal to the distance between $y$
and the graph of $\psi_1$ and
 \begin{align*}
y_d-R+\sqrt{ |\wt y|^2+(R-y_d)^2} &= \frac{|\wt y|^2}
{\sqrt{ |\wt y|^2+(R-y_d)^2} + (R-y_d)}\,\le\, \frac{ |\wt y|^2} {2 (R-y_d)} \le \frac{ |\wt y|^2} {R},
\end{align*}
we have
 \begin{align}
h_{x}(y)-h(y)
\le\int^{y_d}_{R-\sqrt{ |\wt y|^2+(R-y_d)^2}} v(z)dz\le  R^{-1}  |\wt y|^2 \,v(R-\sqrt{ |\wt y|^2+(R-y_d)^2}). \label{e:KG2}
\end{align}
By \eqref{e:KG}-\eqref{e:KG2},
\begin{align*}
I_3
\,\le\,& R^{-1}\int_{E \cap \{  y: y_d \le \delta_D({y})   \}}
 |\wt y|^2 v(y_d)j(|x-y|) dy\\
 &+ R^{-1}\int_{E \cap \{  y:y_d  > \delta_D({y})   \}}
|\wt y|^2  v(R-\sqrt{ |\wt y|^2+(R-y_d)^2}) j(|x-y|) dy=:R^{-1}( L_1+L_2).
\end{align*}

Since
$
E\subset  \{z=(\wt z, z_d)\in \R^d: \ |\widetilde{z}|< R_4
\hbox{ and }  0< z_d\leq 2 R_4\},
$
using polar coordinates for $\wt y$ and by Theorem \ref{t:Jorigin},
\eqref{e:behofV}, \eqref{e:upas2} and \eqref{e:upas3},
we have  that
\begin{eqnarray*}
L_1
&\leq&  c_{8}
\int_{0}^{2 R_4}  v(y_d)\left(
\int_0^{ R_4}
r^2
j((r^2+ |y_d- x_d|^2)^{1/2}) r^{d-2} dr
\right) dy_d\\
&\leq & c_{9}\int_{0}^{2 R_4}  \frac1{(y_d)^{1-\alpha/2}\,
(\ell(y_d^{-2}))^{1/2}}\left(
\int_0^{ R_4} \frac{r^{d}  \ell((r^2+ |y_d- x_d|^2)^{-1}) }
{(r^2+ |y_d- x_d|^2)^{(d+\alpha)/2}} dr
\right) dy_d\\
&\leq & c_{10}\int_{0}^{2 R_4} \frac1{(y_d)^{1-\alpha/2+q} }\left(
\int_0^{ R_4} \frac{r^{d}   }{(r^2+ |y_d- x_d|^2)^{(d+\alpha+2q)/2}}
dr \right) dy_d\\
&\leq & c_{10}\int_{0}^{2 R_4}  \frac1{(y_d)^{1-\alpha/2+q}}\left(
\int_0^{ R_4}
\frac{1}{(r+|y_d-x_d|)^{\alpha+2q}} dr \right) dy_d\\
&\leq& c_{11}\int_{0}^{2 R_4} \frac1{(y_d)^{1-\alpha/2+q}}\left(
\frac1{|y_d-x_d|^{\alpha+2q-1}}
+
\frac{1}{( R_4+|y_d-x_d|)^{\alpha+2q-1}} \right)dy_d
\,\leq\,  c_{12}
\end{eqnarray*}
for some constant $c_8, \dots, c_{12}>0$.
The last inequality
is due to the fact that $q <
(2-\alpha)/20$, which implies $(1-\alpha/2) + \alpha +3q -1 <
(6+7\alpha)/20<1$, so by the dominated convergence theorem,
\bee\label{e:cbdq}
x_d\mapsto  \int_{0}^{2 R_4} \frac1{(y_d)^{1-\alpha/2+q}
|y_d-x_d|^{\alpha+2q-1}} dy_d
\eee
is a strictly positive continuous function in $x_d\in [0,  R_4]$
and hence it is bounded.

On the other hand,  we have, using polar coordinates for $\wt y$,
and by Theorem \ref{t:Jorigin},  \eqref{e:behofV}
and \eqref{e:upas2}--\eqref{e:upas3},
\begin{eqnarray*}
\lefteqn{L_2\leq c_{13}
\int_{0}^{x_d+R_4} \left(
\int_0^{ R_4 \wedge \sqrt{ 2Ry_d-y_d^2}}
 v(R-\sqrt{ r^2+(R-y_d)^2})r^d
j((r^2+ |y_d- x_d|^2)^{1/2}) dr
\right) dy_d}\\
&\leq & c_{14}\int_{0}^{x_d+R_4}\left( \int_0^{ R_4 \wedge \sqrt{
2Ry_d-y_d^2}}  \frac {  (R-\sqrt{ r^2+(R-y_d)^2})^{\alpha/2-1}
\ell((r^2+ |y_d- x_d|^2)^{-1})}{ (\ell((R-\sqrt{
r^2+(R-y_d)^2})^{-2}))^{1/2}(r^2+ |y_d- x_d|^2)^{(d+\alpha)/2}}
r^{d} dr
\right) dy_d\\
&\leq & c_{15}\int_{0}^{x_d+R_4}  \left( \int_0^{R_4 \wedge \sqrt{
2Ry_d-y_d^2}} \frac{r^{d}   }{(R-\sqrt{
r^2+(R-y_d)^2})^{1-\alpha/2+q}\, (r^2+ |y_d-
x_d|^2)^{(d+\alpha+2q)/2}}\ dr
\right) dy_d\\
&\leq &
c_{16}\int_{0}^{x_d+R_4}  \left( \int_0^{R_4 \wedge \sqrt{
2Ry_d-y_d^2}} \frac{1  }{(R-\sqrt{
r^2+(R-y_d)^2})^{1-\alpha/2+q}\, (r+ |y_d- x_d|)^{\alpha+2q}}\ dr
\right) dy_d\, .
\end{eqnarray*}
Since, for $0<r<R_4 \wedge \sqrt{ 2Ry_d-y_d^2}$,
$$
\frac{1  }{R-\sqrt{ r^2+(R-y_d)^2}}
=
\frac{R+\sqrt{ r^2+(R-y_d)^2} }{ (\sqrt{2Ry_d-y_d^2}+r)  (\sqrt{2Ry_d-y_d^2}-r)         }\,\le\,
\frac{c_{17}}{\sqrt{y_d}  (\sqrt{2Ry_d-y_d^2}-r)},
$$
we have
$$
L_2 \le   \int_{0}^{x_d+R_4} \frac{c_{18}  }{(y_d)^{(1-\alpha/2+q)/2}}  \int_0^{ R_4
\wedge \sqrt{ 2Ry_d-y_d^2}} \frac{dr }{ (\sqrt{2Ry_d-y_d^2}-r
)^{1-\alpha/2+q}\, (r+ |y_d- x_d|)^{\alpha+2q}} dy_d\, .
$$

Using the fact that $q \le \frac{\alpha}{20}$,
we see that with $a:=  \sqrt{ 2Ry_d-y_d^2}  $ and $b:=  |y_d- x_d| $,
\begin{align*}
& \int_0^{ R_4 \wedge a} \frac{dr }{ (a-r )^{1-\alpha/2+q}\, (r+
b)^{\alpha+2q}}\\ =& \int_0^{ (R_4 \wedge a)/2} \frac{dr }{ (a-r
)^{1-\alpha/2+q}\, (r+ b)^{\alpha+2q}} +\int_{ (R_4 \wedge a)/2}^{
R_4 \wedge a} \frac{dr }{ (a-r )^{1-\alpha/2+q}\,
(r+ b)^{\alpha+2q}}\\
\le& \frac{2^{1-\alpha/2+q} }{ a^{1-\alpha/2+q}} \int_0^{ (R_4
\wedge a)/2} \frac{dr }{ (r+ b)^{\alpha+2q}} +\frac{1}{(b+ (R_4
\wedge a)/2)^{\alpha+2q}} \int_{ (R_6 \wedge a)/2}^{ R_4 \wedge a}
\frac{dr }{ (a-r )^{1-\alpha/2+q}}\\
\le& \frac{c_{19} }{ a^{1-\alpha/2+q}   b^{(\alpha+2q-1)^+}}
+\frac{c_{19}}{(R_6 \wedge a)^{\alpha+2q}} a^{\alpha/2-q}
\le \frac{c_{20} }{ (y_d)^{(1-\alpha/2+q)/2}   |x_d-y_d|^{(\alpha+2q-1)^+}}
+\frac{c_{20}}{(y_d)^{(\alpha+6q)/4}}\, .
\end{align*}
Thus we obtain
\begin{eqnarray}\label{e:dfsadfs}
L_2
\,\le\, c_{21}\int_{0}^{2R_4}   \frac{dy_d  }{(y_d)^{(1-\alpha/2+q)}
 |y_d- x_d|^{(\alpha+2q-1)^+}
} +c_{21}\int_{0}^{2R_4} \frac{
dy_d}{ (y_d)^{(1+4q)/2}}.
\end{eqnarray}
Since $q<1/10$, the second integral  in \eqref{e:dfsadfs} is bounded. And by the same argument
as the one for \eqref{e:cbdq}, the first integral in \eqref{e:dfsadfs} is also bounded.
We have proved the claim \eqref{e:I1234}. \qed

When $d \ge 2$, define
$\rho_Q (x) := x_d -  \psi_Q (\wt x),$
where $(\wt x, x_d)$ are the coordinates of $x$ in $CS_Q$. Note that
for every $Q \in \partial D$ and $ x \in B(Q, R)\cap D$ we have
\begin{equation}\label{e:d_com}
(1+\Lambda^2)^{-1/2} \,\rho_Q (x) \,\le\, \delta_D(x)  \,\le\,
\rho_Q(x).
\end{equation}

For $a, b>0$, we define
$D_Q( a, b) :=\left\{ y\in D: a >\rho_Q(y) >0,\, |\wt y | < b
\right\}$ when $d \ge2$. When $d=1$, we simply take $D_Q( a, b)=D_Q( a):=B(Q,a) \cap D$.

\begin{lemma}\label{L:2}
There are constants
$R_5=R_5( R, \Lambda, \alpha, \ell)\in (0, R_4/(4 \sqrt{1+(1+ \Lambda)^2}))$ and
$C_i=C_i(R,  \Lambda, \alpha)>0$, $i=17, 18$, such
that for every $r \le R_5$, $Q \in
\partial D$ and $x \in D_Q ( r, r)$,
 \bee\label{e:L:2}
\P_{x}\left(X_{ \tau_{ D_Q ( r, r)}} \in  D\right) \ge {C_{17}} V( \delta_D (x))
 \eee
and
 \bee\label{e:L:3}
\E_x\left[ \tau_{ D_Q  ( r, r )}\right]\,\le\, {C_{18}} V(\delta_D (x)).
 \eee
\end{lemma}

\pf
Without loss of generality, we assume $Q=0$.
For $d \ge2$,  let
$\psi=\psi_0:\bR^{d-1}\to \bR$ be the $C^{1,1}$ function  and $CS_0$ be the coordinate system
in the definition of $C^{1, 1}$ open set so that
$B(0, R)\cap D= \big\{(\wt y, \, y_d) \in B(0, R)\textrm{ in } CS_0: y_d >
\psi (\wt y) \big\}.$
Let $\rho (y) := y_d -  \psi (\wt y)$ and $D ( a, b):=D_0 ( a, b)$ for $d \ge2$.
When $d=1$, $D ( a, b)$ is simply $B(0, a) \cap D$. The remainder of the proof is written for $d\ge 2$. The interpretation in the case $d=1$ is obvious.

Note that
\bee\label{e:lsd}
|y|^2 = |\wt y|^2+ |y_d|^2 <r^2 +(|y_d- \psi(\wt y)|+ |\psi(\wt
y)|)^2 < (1+(1+ \Lambda)^2) r^2 \quad \text{for every } y \in
D(r,r)\, .
\eee
Hence, by letting $\wh r :=R_4/\sqrt{1+(1+ \Lambda)^2}$,
$D ( r, s) \subset D(\wh r, \wh r)\subset B(0, R_4)
\cap D \subset B(0, R)\cap D$  for every $ r,s \le \wh r .
$
Define
$$
h(y) := V(\delta_D(y))
{\bf 1}_{B(0, R) \cap D}(y).$$

Let $g$ be a non-negative smooth radial function with compact support
such that $g(x)=0$ for
$|x|>1$ and $\int_{\R^d} g (x) dx=1$.
For $k\geq 1$, define $g_k(x)=2^{kd} g (2^k x)$ and
$$
h^{(k)}(z):= ( g_k*h)(z) :=\int_{\R^d}g_k (y) h(z-y)dy\, ,
$$
and let $B_k:=\left\{x \in D \cap B(0, R_4): \delta_{D \cap B(0, R_4)    }(x) \ge 2^{-k}\right\}$.
Since $h^{(k)}$ is $C^{\infty}$, $\AA h^{(k)}$ is well defined everywhere. We claim that
\begin{equation} \label{e:*333}
-C_{16} \le \AA h^{(k)} \leq  C_{16} \quad \text{ on } B_k\, ,
\end{equation}
where $C_{16}$ is the constant from Lemma \ref{L:Main}. Indeed, for $x\in B_k$ and $z\in B(0,2^{-k})$ it holds that $x-z\in D \cap B(0, R_4)$. Hence, by Lemma \ref{L:Main} the following limit exists:
\begin{eqnarray*}
\lefteqn{\lim_{\varepsilon\to 0} \int_{|x-y|>\varepsilon} \left(h(y-z)-h(x-z)\right)\, j(|x-y|)\, dy}\\
&=&\lim_{\varepsilon\to 0} \int_{|(x-z)-y'|>\varepsilon} \left(h(y')-h(x-z)\right)\, j(|(x-z)-y'|)\, dy'
\,=\,\AA h(x-z)\, .
\end{eqnarray*}
Moreover, by the same Lemma \ref{L:Main} it holds that $-C_{16}\le \AA h(x-z) \le C_{16}$. Next,
\begin{eqnarray*}
\lefteqn{\int_{|x-y|>\varepsilon} (h^{(k)}(y)-h^{(k)}(x))\, j(|x-y|)\, dy}\\
&=&\int_{|x-y|>\varepsilon} \left(\int_{\R^d}g_k(z)(h(y-z)-h(x-z))\, dz\right)\, j(|x-y|)\, dy\\
&=&\int_{|z|<2^{-k}}g_k(z) \left(\int_{|x-y|>\varepsilon}\left(h(y-z)-h(x-z)\right)\, j(|x-y|)\, dy\right)\, dz.
\end{eqnarray*}
By letting $\varepsilon \to 0$ and using the dominated convergence theorem, it follows that
$$
\AA h^{(k)}(x)=\int_{|z|<2^{-k}} g_k(z)\AA h(x-z)\, dz\le  C_{16}\int_{|z|<2^{-k}} g_k(z) \, dz =C_{16}\, .
$$
The left-hand side inequality in \eqref{e:*333} is obtained in the same way.

Using the fact that $\AA$ restricted to $C^\infty_c$
coincides with the infinitesimal generator of the process $X$, we see that
the following Dynkin formula is true;
for $f \in C_c^{\infty}(\R^d)$ and any bounded open subset $U$ of $\R^d$,
\begin{equation} \label{e:*334}
\E_x\int_0^{\tau_U}  \AA f(X_t) dt
=\E_x[f(X_{\tau_U})]- f(x).
\end{equation}

Let
$U\subset  D \cap B(0, R_4)$. By using \eqref{e:*334} for $U\cap
B_k$ and $h^{(k)}$, the estimates \eqref{e:*333}, the fact that $h^{(k)}$ are in
$C^\infty_c(\R^d)$, and by letting $k\to \infty$ we get
 \bee\label{e:sh4} h(x) \ge
\E_{x}\left[h(X_{\tau_U})\right]-C_{16}\E_x[\tau_U] \quad \text{and}
\quad h(x) \le \E_{x}\left[h(X_{\tau_U})\right]+C_{16} \E_x[\tau_U].
 \eee
Now, we have by \eqref{e:d_com} and \eqref{e:sh4}, for every
$\lambda \ge 1$ and $x \in D(\lambda^{-1}  \wh r , \lambda^{-1}  \wh
r )$,
\begin{align}
&V(\delta_D(x)) \,=\, h(x) \nn\\
&\ge \E_{x}\left[h\left(X_{ \tau_{D(\lambda^{-1}  \wh r ,
\lambda^{-1}  \wh r )}} \right); X_{ \tau_{D(\lambda^{-1}  \wh r ,
\lambda^{-1}  \wh r )}} \in D(  \wh r , \lambda^{-1}  \wh r )
\setminus
D(\lambda^{-1}  \wh r , \lambda^{-1}  \wh r ) \right] - C_{16}
\E_x\left[\tau_{D(\lambda^{-1}  \wh r , \lambda^{-1}  \wh r )}\right]\nn\\
&\ge
  V\big(\lambda^{-1}(1+\Lambda^2)^{-1/2} \wh r \big)\P_{x}
  \left(X_{ \tau_{D(\lambda^{-1}  \wh r , \lambda^{-1}
  \wh r )}} \in D(  \wh r , \lambda^{-1}  \wh r ) \setminus
D(\lambda^{-1}  \wh r , \lambda^{-1}  \wh r ) \right)-C_{16}
\E_x\left[\tau_{D(\lambda^{-1}  \wh r , \lambda^{-1}  \wh r )}\right] \label{e:ss}.
\end{align}
We also have from \eqref{e:sh4}
\begin{eqnarray}
V(\delta_D(x)) =
h(x) &\le &\E_{x}\left[h \left( X_{  \tau_{D(\lambda^{-1}  \wh r ,
\lambda^{-1}  \wh r )}} \right)
\right]+ C_{16} \E_x[\tau_{D(\lambda^{-1}  \wh r , \lambda^{-1}
\wh r )}]\nonumber \\
&\le &
V(R) \P_{x}\left(X_{ \tau_{D(\lambda^{-1}  \wh r ,
\lambda^{-1}  \wh r )}} \in D \right) + C_{16} \E_x[\tau_{D(\lambda^{-1}
\wh r , \lambda^{-1}  \wh r )}] \label{e:ss1}.
\end{eqnarray}
By
\eqref{newls} and the monotonicity of $j$, for every $\lambda
\ge 4$ and $x \in D(\lambda^{-1}  \wh r , \lambda^{-1}  \wh r )$,
\begin{eqnarray*}
\P_{x}\left(X_{ \tau_{D(\lambda^{-1}  \wh r , \lambda^{-1}  \wh r )}}
\in D \right)
&\ge & \P_{x}\left(X_{ \tau_{D(\lambda^{-1}  \wh r , \lambda^{-1}
\wh r )}} \in D(  \wh r , \lambda^{-1} \wh r ) \setminus
D(\lambda^{-1}  \wh r , \lambda^{-1}  \wh r )\right) \\
&\ge &
\P_{x}\left(X_{ \tau_{D(\lambda^{-1}  \wh r ,\lambda^{-1} \wh r )}}
\in D( 3\lambda^{-1} \wh r ,\lambda^{-1} \wh r ) \setminus
D(2\lambda^{-1}  \wh r , \lambda^{-1}  \wh r ) \right)\\
&=&
\E_x \left[ \int_0^{ \tau_{D(
\lambda^{-1}  \wh r , \lambda^{-1}  \wh r )}} \int_{D(
3\lambda^{-1} \wh r ,\lambda^{-1} \wh r ) \setminus
D(2\lambda^{-1}  \wh r , \lambda^{-1}  \wh r )}j(|X_s-y|) dy ds\right]\\
&\ge &
\left(\int_{D( 3\lambda^{-1} \wh r ,\lambda^{-1} \wh r ) \setminus
D(2\lambda^{-1}  \wh r , \lambda^{-1}  \wh r )} dy \right)
j(10\lambda^{-1}  \wh r ) \, \E_x \left[
\tau_{D(\lambda^{-1}  \wh r , \lambda^{-1}  \wh r )}\right]\\
&\ge &
c_1 (\lambda^{-1} \wh r )^d j(10\lambda^{-1}  \wh r ) \, \E_x \left[
\tau_{D(\lambda^{-1}  \wh r , \lambda^{-1}  \wh r )}\right].
\end{eqnarray*}
Now, applying Theorem \ref{t:Jorigin}, we get for  $x \in
D(\lambda^{-1}  \wh r , \lambda^{-1}  \wh r )$
\begin{align}
\P_{x}\left(X_{ \tau_{D(\lambda^{-1}  \wh r , \lambda^{-1}  \wh r )}} \in D \right) \ge c_2
{\ell \big((10\lambda^{-1}  \wh r )^{-2}\big)}{\lambda^\alpha}\E_x \left[
\tau_{D(\lambda^{-1}  \wh r , \lambda^{-1}  \wh r )}\right].\label{e:ss2}
\end{align}
Thus from \eqref{e:ss}--\eqref{e:ss2}, for every $x \in
D(\lambda^{-1} \wh r , \lambda^{-1}  \wh r )$
\begin{align}
V(\delta_D(x))\ge \left(c_2V\big(\lambda^{-1}(1+\Lambda^2)^{-1/2} \wh r \big) {\ell
\big((10\lambda^{-1}  \wh r )^{-2}\big)}{\lambda^\alpha}    -C_{16}\right)
\E_x[\tau_{D(\lambda^{-1}  \wh r , \lambda^{-1}  \wh r )}]
\label{e:ss4}
\end{align}
and
\begin{align}
V(\delta_D(x)) \le
c_3\left(1+ \left(\ell \big((10\lambda^{-1}  \wh r )^{-2}\big)\right)^{-1}\lambda^{-\alpha} \right)
\P_{x}\left(X_{
\tau_{D(\lambda^{-1}  \wh r , \lambda^{-1}  \wh r )}} \in D \right)
\label{e:ss5}.
\end{align}

Using first \eqref{e:behofV} and then Potter's
Theorem (\cite[Theorem 1.5.6 (1)]{BGT}), we see that there exists a
large $\lambda_0 >4$ such that  for every $\lambda \ge \lambda_0$
\begin{align}
& V\big(\lambda^{-1}(1+\Lambda^2)^{-1/2}\wh r \big) {\ell \big((10\lambda^{-1}
\wh r )^{-2}\big)}{\lambda^\alpha}\nn \\&
\ge
c_4 \wh r ^{\alpha/2}
    (1+\Lambda^2)^{-\alpha/4}  \lambda^{{\alpha/2}}
     \left(\ell\Big(\big(\lambda^{-1}(1+\Lambda^2)^{-1/2}\wh r
     \big)^{-2}\Big)\right)^{-1/2}
    {\ell \big((10\lambda^{-1}  \wh r )^{-2}\big)} \ge 2C_{16}/c_2\, \label{e:ss6}
\end{align}
and
\begin{align}\label{e:ss7}
\left(\ell \big((10\lambda^{-1}  \wh r )^{-2}\big)\right)^{-1}\lambda^{-\alpha}  \le c_5.
\end{align}
Combining
\eqref{e:ss4}--\eqref{e:ss7}, we have proved the lemma
with $R_5:=\lambda_0^{-1} \wh r$. \qed

It is clear that every $C^{1,1}$ open set is $\kappa$-fat, i.e.,
for any $C^{1, 1}$ open set with $C^{1,1}$ characteristics
$(R, \Lambda)$, there exists a  constant $\kappa\in (0,1/2]$,
which depends only on
$(R, \Lambda)$, such that for each $Q \in \partial D$
and $r \in (0, R)$, $D \cap B(Q,r)$
contains a ball $B(A_r(Q),\kappa r)$
of radius $\kappa r$.
In the rest of this paper, whenever we deal with
$C^{1,1}$ open sets, the constants $\Lambda$, $R$ and $\kappa$ will
have the meaning described above.

Recall that $g$ is defined in \eqref{d:gz0}.

\begin{thm}\label{t:Ge} Suppose that
$D$ is a bounded $C^{1,1}$ open set in $\R^d$ with $C^{1,1}$ characteristics
$(R, \Lambda)$. Then there exists
$C_{19}=C_{19}(R, \Lambda, \alpha,\textrm{diam}(D))>0$ such that
\begin{equation}\label{e:z1}
C_{19}^{-1} \,\left(V(\delta_D(x))
   \wedge 1 \right)\,\le \,g(x)\, \le\,
C_{19}\,\,\left( V(\delta_D(x))   \wedge 1 \right), \qquad x \in D,
\end{equation}
or equivalently there exists $C_{20}=C_{20}(R, \Lambda, \alpha, \textrm{diam}(D))>0$ such that
\begin{equation}\label{e:z2}
C_{20}^{-1} \,\left( \frac{(\delta_D(x))^{\alpha/2}}
{\sqrt{\ell((\delta_D(x)  )^{-2})}}      \wedge 1 \right)\,\le \,g(x)\, \le\,
C_{20}\,\,\left( \frac{(\delta_D(x))^{\alpha/2}}
{\sqrt{\ell((\delta_D(x)  )^{-2})}}      \wedge 1 \right), \qquad x \in D.
\end{equation}

\end{thm}

\pf
Since $d=1$ case is simpler, we give the proof for $d \ge 2$ only.
Recall that $R_3$ is the constant in \eqref{e:rsmall} and
$\eps_1= R_3\kappa /24$. Since
$g(x)= G_D(x, z_0) \wedge C_{13}$ and
$g(x)= G_D(x, z_0)$ for $\delta_D(x) < 6 \eps_1$, it suffices to
show that there exist  $r^* \in (0, 6 \eps_1)$ and $c_1>1$ such that
 \bee \label{e:fclaim1} c_1^{-1}  \frac{(\delta_D(x))^{\alpha/2}}
{\sqrt{\ell((\delta_D(x)  )^{-2})}}     \,\le \,
G_D(x, z_0)
\, \le\, c_1\,
\frac{(\delta_D(x))^{\alpha/2}} {\sqrt{\ell((\delta_D(x)  )^{-2})}}
, \qquad  \delta_D(x) < r^*\, .
 \eee

Let
$r^*:= (R_5/4) \wedge (\eps_1/(4\sqrt{1+(1+ \Lambda)^2}))$
 and suppose that $\delta_D(x) < r^*$.
Choose $x_0\in\partial D$ satisfying $\delta_D(x)=|x-x_0|$. We fix
the $C^{1, 1}$ function $\psi$ and the coordinate system $CS=CS_{x_0}$ in the definition of $C^{1, 1}$ open set
so that $x=( \wt 0, x_d)$ with
$0<x_d <r^*$,
$$
B(x_0, R)\cap D=
\{ y=(\wt y, \, y_d)\in B(0, R) \textrm{ in } CS : y_d > \psi (\wt y) \}.
$$

Let $x_1:=(\wt 0, r^*/2)$ and $D_*:=D(r^*,r^*)=\{ y\in D: r^* > y_d -  \psi (\wt y) >0,\, |\wt y | < r^*\}$.
Since $B(x_1, c_2 r^*)  \subset D_*$ for small $c_2>0$,
by Theorem \ref{t:Gorigin}, Theorem \ref{BHP} and the fact that $D$ is bounded,
\begin{eqnarray*}
G_D(x, z_0) &\le& c_3 G_D(x_1, z_0)
\frac{\P_{x}\left(X_{ \tau_{ D_*}} \in  B(z_0, \eps_1/4)\right)}
{\P_{x_1}\left(X_{ \tau_{ D_*}} \in  B(z_0, \eps_1/4)\right)} \\
&\le& c_3 \,G(x_1, z_0)
\frac{\P_{x}\left(X_{ \tau_{ D_*}} \in  B(z_0, \eps_1/4)\right)}
{\P_{x_1}\left(X_{ \tau_{ B(x_1, c_2 r^*) }} \in  B(z_0, \eps_1/4)\right)}
\\&\le&
c_4 \,\P_{x}\left(X_{ \tau_{ D_*}} \in  B(z_0, \eps_1/4)\right)\,
\le\,  c_5\, \E_x\left[ \tau_{ D_*}\right]
\end{eqnarray*}
where in the last inequality we used \eqref{newls}
and the fact that dist$(D_*, B(z_0, \eps_1/4)) \ge \delta_D(z_0) - \eps_1/4 -\sqrt{1+(1+ \Lambda)^2} r^* \ge \eps_1$  (see \eqref{e:lsd}). On the other hand, by Theorem \ref{BHP},  Lemma \ref{G:g3} and the fact that $D$ is bounded,
\begin{eqnarray*}
G_D(x, z_0) \ge c_5 G_D(x_1, z_0)\frac{ \P_{x}\left(X_{ \tau_{ D_*}}
\in  D\right)}{\P_{x_1}\left(X_{ \tau_{ D_*}} \in   D\right)}\ge c_6 \P_{x}\left(X_{ \tau_{ D_*}}   \in  D\right).
\end{eqnarray*}
By applying \eqref{e:L:2}--\eqref{e:L:3}, we have proved
\eqref{e:z1}. The inequalities \eqref{e:z2} follow from \eqref{e:behofV}. \qed

Now we give the proof of Theorem \ref{t:Gest2}, which is the main result of this paper.

\medskip
\noindent {\bf Proof of  Theorem \ref{t:Gest2}}.
By \cite[Theorem 1.5.3]{BGT}, the local boundedness and strict
positivity of $\ell$, there exists $c_1>1$ such that for every
$0<s\le t \le 6\textrm{ diam}(D)=:6M$
$$
\frac{s^\alpha}{(\ell(s^{-2}))} \,\le\,
\frac{c_1t^\alpha}{(\ell(t^{-2}))}.
$$ Thus, if $s,t,u \le 6M$ then
 \beq\label{e:sim0} c_1^{-1}
\left(\frac{s^\alpha}{\ell(s^{-2})} \vee
\frac{t^\alpha}{\ell(t^{-2})} \vee \frac{u^\alpha}{\ell(u^{-2})}
\right)\,\le\, \frac{(s\vee t \vee u)^\alpha}{\ell((s\vee t \vee
u)^{-2})} \,\le\,
c_1\left(\frac{s^\alpha}{\ell(s^{-2})} \vee
\frac{t^\alpha}{\ell(t^{-2})} \vee \frac{u^\alpha}{\ell(u^{-2})}
\right).
 \eeq

Also note that
 by \eqref{d:gz1}--\eqref{e:areq}, for every $A\in\BB(x, y)$
 \beq\label{e:Gest2}
\frac12 \wedge \frac{\min_{\kappa R_5 \le s \le R_5} \ell(s^{-2})}{\max_{\eps_1 \le s \le M} \ell(s^{-2})} \, \le \,
\frac{\ell((\delta_D(A))^{-2})} {\ell((\delta_D(x) \vee \delta_D (y)
\vee |x-y|)^{-2})} \,\le\, 2 \vee \frac{\max_{\kappa R_5 \le s \le R_5} \ell(s^{-2})}{\min_{\eps_1 \le s \le M}
\ell(s^{-2})}
 \eeq
 and
 \beq\label{e:Gee1} \left(\frac16\wedge
\frac{\eps_1}{R_3}\right)\delta_D(A) \,\le\,\delta_D(x) \vee
\delta_D (y) \vee |x-y|\, \le\, 2\kappa^{-1}\left(\frac{M}{2R_3}
\vee 1\right) \delta_D(A).
 \eeq

By Theorem \ref{t:Gest} and Theorem \ref{t:Ge}, we have that
\begin{eqnarray*}
&&c_2^{-1}\,\frac{(\delta_D(x) \delta_D(y))^{\alpha/2}
\ell((\delta_D(A))^{-2})  }{\sqrt{\ell((\delta_D(x)  )^{-2})
\ell((\delta_D(y)  )^{-2})}(\delta_D(A))^{\alpha}}
\frac1{\ell(|x-y|^{-2})|x-y|^{d-\alpha}}
\,\le\, G_D(x,y)
\\
&\le&
c_2\,\,\frac{(\delta_D(x) \delta_D(y))^{\alpha/2}
\ell((\delta_D(A))^{-2})  }{\sqrt{\ell((\delta_D(x)
)^{-2})\ell((\delta_D(y)  )^{-2})}(\delta_D(A))^{\alpha}}
\frac1{\ell(|x-y|^{-2})|x-y|^{d-\alpha}}.
\end{eqnarray*}
Since \eqref{e:sim0}--\eqref{e:Gee1} and the identity
$$\frac{a}{b} \wedge \frac{b}{a} \wedge\frac{ab}{c^2}= \frac{ab}{(a \vee b \vee c)^2}, \quad \text{ for all } a,b,c>0$$ imply that
\begin{align*}
&c_3^{-1}\frac{(\delta_D(x) )^{\alpha/2} \sqrt{\ell((\delta_D(y)
)^{-2})}}{(\delta_D(y) )^{\alpha/2}\sqrt{\ell((\delta_D(x)
)^{-2})}} \wedge \frac{(\delta_D(y) )^{\alpha/2}
\sqrt{\ell((\delta_D(x)  )^{-2})}} {(\delta_D(x) )^{\alpha/2}
\sqrt{\ell((\delta_D(y)  )^{-2})}  } \wedge \frac{(\delta_D(x)
\delta_D(y))^{\alpha/2}  \ell(|x-y|^{-2})  }{\sqrt{\ell(
(\delta_D(x) )^{-2})\ell((\delta_D(y)  )^{-2})}|x-y|^{\alpha}}
\\
\le&
\frac{(\delta_D(x) \delta_D(y))^{\alpha/2}
\ell((\delta_D(A))^{-2})  }{\sqrt{\ell((\delta_D(x)  )^{-2})
\ell((\delta_D(y)  )^{-2})}(\delta_D(A))^{\alpha}}
\\
\le&
c_3\frac{(\delta_D(x) )^{\alpha/2} \sqrt{\ell((\delta_D(y)
)^{-2})}}{(\delta_D(y) )^{\alpha/2}\sqrt{\ell((\delta_D(x)
)^{-2})}} \wedge \frac{(\delta_D(y) )^{\alpha/2}
\sqrt{\ell((\delta_D(x)  )^{-2})}} {(\delta_D(x) )^{\alpha/2}
\sqrt{\ell((\delta_D(y)  )^{-2})}  } \wedge \frac{(\delta_D(x)
\delta_D(y))^{\alpha/2}  \ell(|x-y|^{-2})  }{\sqrt{\ell(
(\delta_D(x) )^{-2})\ell((\delta_D(y)  )^{-2})}|x-y|^{\alpha}},
\end{align*}
to prove the theorem, we only need to show
\begin{align}
&\frac{(\delta_D(x) )^{\alpha/2} \sqrt{\ell((\delta_D(y)
)^{-2})}}{(\delta_D(y) )^{\alpha/2}\sqrt{\ell((\delta_D(x)
)^{-2})}} \wedge \frac{(\delta_D(y) )^{\alpha/2}
\sqrt{\ell((\delta_D(x)  )^{-2})}} {(\delta_D(x) )^{\alpha/2}
\sqrt{\ell((\delta_D(y)  )^{-2})}  } \wedge \frac{(\delta_D(x)
\delta_D(y))^{\alpha/2}  \ell(|x-y|^{-2})  }{\sqrt{\ell(
(\delta_D(x) )^{-2})\ell((\delta_D(y)  )^{-2})}|x-y|^{\alpha}}\nonumber\\
\le&\left(1 \wedge\frac{(\delta_D(x) \delta_D(y))^{\alpha/2}
\ell(|x-y|^{-2})  }{\sqrt{\ell((\delta_D(x)  )^{-2})
\ell((\delta_D(y)  )^{-2})}|x-y|^{\alpha}} \right)\nonumber\\
\le& c_5\,\,\left(\frac{(\delta_D(x) )^{\alpha/2}
\sqrt{\ell((\delta_D(y) )^{-2})}  }{(\delta_D(y)
)^{\alpha/2}\sqrt{\ell((\delta_D(x)  )^{-2})}} \wedge
\frac{(\delta_D(y) )^{\alpha/2}\sqrt{\ell((\delta_D(x)  )^{-2})}}
{(\delta_D(x) )^{\alpha/2} \sqrt{\ell((\delta_D(y)  )^{-2})}  }
\wedge \frac{(\delta_D(x) \delta_D(y))^{\alpha/2}  \ell(|x-y|^{-2})
}{\sqrt{\ell((\delta_D(x)  )^{-2})\ell((\delta_D(y)
)^{-2})}|x-y|^{\alpha}}\right).\nonumber
\\ \label{e:sim1}
\end{align}
Since the first inequality is clear, we will proceed to the second inequality.

By symmetry, we only need to consider the cases $\delta_D(y)\le \delta_D(x)/3$ and $\delta_D(x)/3 \le \delta_D(y)\le 3 \delta_D(x)$, and, using the fact $\ell$ is slowly varying and \cite[Theorem 1.5.3]{BGT}, the case $\delta_D(x)/3 \le \delta_D(y)\le 3 \delta_D(x)$ is clear.

Now we assume $\delta_D(y)\le \delta_D(x)/3$. Since $|x-y| \ge \delta_D(x)-\delta_D(y)  \ge 2 \delta_D(x)/3$ in this case, using
\cite[Theorem 1.5.3]{BGT}, the continuity and strict
positivity of $\ell$, \begin{align*}
& \frac{\delta_D(x) \delta_D(y)  \ell(|x-y|^{-2})^{2/\alpha}  }{(\ell((\delta_D(x)  )^{-2}))^{1/\alpha}(\ell((\delta_D(y)  )^{-2}))^{1/\alpha}|x-y|^{2}}
\\
&\le c_6 \frac{\delta_D(x) \delta_D(y)   }{(\ell((\delta_D(x)  )^{-2}))^{1/\alpha}(\ell((\delta_D(y)  )^{-2}))^{1/\alpha}}
\left(\frac{(\ell((\delta_D(x)  )^{-2}))^{1/\alpha}}{\delta_D(x)}\right)^{2}
=c_6 \frac{ \delta_D(y)   }{(\ell((\delta_D(y)  )^{-2}))^{1/\alpha}}
\frac{(\ell((\delta_D(x)  )^{-2}))^{1/\alpha}}{\delta_D(x)}
 \end{align*}
 and
$$
\frac{(\delta_D(y) )^{\alpha/2}\sqrt{\ell((\delta_D(x)  )^{-2})}} {(\delta_D(x) )^{\alpha/2} \sqrt{\ell((\delta_D(y)  )^{-2})}  }\le c_7\frac{(\delta_D(x) )^{\alpha/2}\sqrt{\ell((\delta_D(y)  )^{-2})}} {(\delta_D(y) )^{\alpha/2} \sqrt{\ell((\delta_D(x)  )^{-2})}  }\, .
$$
Thus
\begin{eqnarray*}
&&\left(1 \wedge\frac{(\delta_D(x) \delta_D(y))^{\alpha/2}  \ell(|x-y|^{-2})  }{\sqrt{\ell((\delta_D(x)  )^{-2})\ell((\delta_D(y)  )^{-2})}|x-y|^{\alpha}} \right)\\
 &&\le c_8\left(\frac{(\delta_D(y) )^{\alpha/2}\sqrt{\ell((\delta_D(x)  )^{-2})}} {(\delta_D(x) )^{\alpha/2} \sqrt{\ell((\delta_D(y)  )^{-2})}}
 \wedge \frac{(\delta_D(x) \delta_D(y))^{\alpha/2}  \ell(|x-y|^{-2})  }{\sqrt{\ell((\delta_D(x)  )^{-2})\ell((\delta_D(y)  )^{-2})}|x-y|^{\alpha}}\right)\\
 &&\le c_9\left(\frac{(\delta_D(x) )^{\alpha/2}
\sqrt{\ell((\delta_D(y) )^{-2})}  }{(\delta_D(y)
)^{\alpha/2}\sqrt{\ell((\delta_D(x)  )^{-2})}} \wedge
\frac{(\delta_D(y) )^{\alpha/2}\sqrt{\ell((\delta_D(x)  )^{-2})}}
{(\delta_D(x) )^{\alpha/2} \sqrt{\ell((\delta_D(y)  )^{-2})}  }
\wedge \frac{(\delta_D(x) \delta_D(y))^{\alpha/2}  \ell(|x-y|^{-2})
}{\sqrt{\ell((\delta_D(x)  )^{-2})\ell((\delta_D(y)
)^{-2})}|x-y|^{\alpha}}\right).
\end{eqnarray*}
We have obtained the second inequality in \eqref{e:sim1}.
\qed

\begin{remark}\label{r:alternative form}
Similarly as in Theorem \ref{t:Ge}, by use of \eqref{e:behofV},
the inequalities \eqref{e:Gest21} imply the alternative forms given in
\eqref{e:Gest21-alt1} and \eqref{e:Gest21-alt2}.
\end{remark}

Now we give the proof of Theorem \ref{t:bhp}, which is a consequence of Theorems
\ref{t:Gest2} and
\ref{BHP}.

\medskip

\noindent {\bf Proof of  Theorem \ref{t:bhp}}.
Using the interior ball condition of $D$,
the following holds: For every $Q\in \partial D$ and $r\le R$
there is a ball $B=B(z^r_Q, r)$ of radius $r$ such that $B\subset
 D$ and $\partial B \cap \partial D=\{Q\}$.
In addition, it follows from \cite[Lemma 2.2]{S} that, for each $Q \in
\partial D$, we can choose a constant $c_2=c_2(d, \Lambda) \in
(0,1/8]$ and a bounded $C^{1,1}$ open set $U_Q$ with uniform
characteristics $(R_*, \Lambda_*)$ depending
only on $(R, \Lambda)$ and $d$
such that $B(Q, c_2 R) \cap D \subset U_Q \subset B(Q, R) \cap
D$ and
\begin{equation}\label{e:dcp}
\delta_D(y)=\delta_{U_Q}(y) \quad \text{ for every }y \in B(Q, c_2 R) \cap D.
\end{equation}

Assume that $r \in (0, c_2R]$, $Q\in \partial D$ and $u$ is
nonnegative function in $\R^d$ harmonic in
$D \cap B(Q,
r)=U_Q \cap B(Q,
r)$ with respect to $X$ and vanishes continuously on $ D^c \cap
B(Q, r)$. Let $z_Q:=z_Q^{c_2R}$. By \cite[Lemma 4.2]{CKSV} and its proof,
we see that $u$ and $x \to G_{U_Q}(x, z_Q)$ are regular harmonic in
$U_Q \cap B(Q,
2r/3)$ with respect to $X$. Since the $C^{1,1}$ characteristics of $U_Q$ depend
only on $(R, \Lambda)$ and $d$,
by the boundary Harnack principle (Theorem \ref{BHP}),
there exist
$r_1=r_1(\alpha, \ell, R, \Lambda) \in (0, 1/4]$ and $c_3=c_3(\alpha, \ell, R, \Lambda)>0$
such that for any $r\in (0, r_1]$ we have
$$
\frac{u(x)}{u(y)} \le c_3 \frac{G_{U_Q}(x, z_Q)}{G_{U_Q}(y,
z_Q)}
\quad \text{ for every } x,y \in B(Q, r/2) \cap D.
$$

Now applying Theorem \ref{t:Gest2} to $G_{U_Q}(x, z_Q)$ and
$G_{U_Q}(y, z_Q)$, then using \eqref{e:dcp},
we conclude that for $ r \in (0, (c_2R \wedge r_1)]$
$$
\frac{u(x)}{u(y)} \le c_4
\frac{\delta^{\alpha/2}_{U_Q}(x)\sqrt{\ell((\delta_{U_Q}(y))^{-2})}}{\delta^{\alpha/2}_{U_Q}(y)
\sqrt{\ell((\delta_{U_Q}(x))^{-2})}}= c_4
\frac{\delta^{\alpha/2}_D(x)\sqrt{\ell((\delta_{D}(y))^{-2})}}{\delta^{\alpha/2}_D(y)
\sqrt{\ell((\delta_{D}(x))^{-2})}}\quad \text{
for every } x,y \in B(Q, r/2) \cap D
$$
for some $c_4=c_4(\alpha, \ell, R, \Lambda)>0$.
The form \eqref{e:bhp_m} given in the statement of the theorem is
equivalent to the one in the display above for
$ r \in (0, (c_2R \wedge r_1)]$. Now the case $ r \in ((c_2R \wedge r_1),
(R \wedge 1)/4]$ follows from the case $ r \in (0, (c_2R \wedge r_1)]$
and Theorem \ref{HP2}.
\qed

\vspace{.1in}
\begin{singlespace}

\end{singlespace}
\end{doublespace}

\bigskip

{\bf Panki Kim}

Department of Mathematical Sciences and Research Institute of Mathematics,

Seoul National University,
San56-1 Shinrim-dong Kwanak-gu,
Seoul 151-747, Republic of Korea

E-mail: \texttt{pkim@snu.ac.kr}

\bigskip

{\bf Renming Song}

Department of Mathematics, University of Illinois, Urbana, IL 61801, USA

E-mail: \texttt{rsong@math.uiuc.edu}

\bigskip

{\bf Zoran Vondra{\v{c}}ek}

Department of Mathematics, University of Zagreb, Zagreb, Croatia

Email: \texttt{vondra@math.hr}

\end{document}